\documentclass[12pt]{article}
\usepackage{amsmath,amssymb,a4wide,amsthm}
\usepackage{hyperref}
\usepackage{enumitem}
\usepackage{xcolor}
\usepackage{comment}

\newtheorem{proposition}{Proposition}
\newtheorem{corollary}{Corollary}
\newtheorem{lemma}{Lemma}
\newtheorem{remark}{Remark}
\newtheorem{discussion}{Discussion}
\newtheorem{example}{Example}
\newtheorem{assumption}{Assumption}
\usepackage{graphicx}
\usepackage{hyperref}

\allowdisplaybreaks

\def\R{\mathbb{R}}
\def\N{\mathbb{N}}
\def\Cc{\mathcal{C}}
\def\Vc{\mathcal{V}}

\def\Yc{\mathcal{Y}}
\def\Wc{\mathcal{W}}
\def\Nc{\mathcal{N}}
\def\Lc{\mathcal{L}}
\def\zb{\bar{\zeta}}
\def\ub{\bar{u}}
\def\fb{\bar{f}}

\def\ud{{u^\dagger}}
\def\fd{{f^\dagger}}
\def\lomy{{L^\infty(\Omega_{y^\dagger})}}
\def\too{-\text{to}-}
\def\tom{(0,T)\times\Omega}
\def\source{\varphi}
\def\srcWave{\widetilde{\varphi}}
\def\indmeas{m} 
\def\partostate{P} 
\def\spaceA{A}
\def\spaceB{B}
\newcommand{\compt}{\hookrightarrow\mathrel{\mspace{-15mu}}\rightarrow
}
\newcommand{\embed}{\hookrightarrow}
\newcommand{\blue}[1]{\textcolor{black}{#1}}

\title{Discretization of parameter identification in PDEs using Neural Networks}
\author{Barbara Kaltenbacher, University of Klagenfurt \\and\\ Tram Thi Ngoc Nguyen, University of Graz}
\date{ }

\begin{document}

\maketitle

\paragraph{Abstract.} We consider the ill-posed inverse problem of \blue{identifying a nonlinearity in a time-dependent PDE model. The nonlinearity is approximated by a neural network, and needs to be determined alongside other unknown physical parameters and the unknown state. Hence, it is not possible to construct input-output data pairs to perform a supervised training process. Proposing an all-at-once approach, we bypass the need for training data and recover all the unknowns simultaneously.} In the general case, the approximation via a neural network can be realized as a discretization scheme, \blue{and the training with noisy data can be viewed as an ill-posed inverse problem}. Therefore, we study discretization of regularization in terms of Tikhonov and projected Landweber methods for \blue{discretization of} inverse problems, and prove convergence when the discretization error \blue{(network approximation error)} and the noise level tend to zero.\\[1ex]

\section{Introduction}\label{sec:intro}

Parameter identification in partial differential equations (PDEs) from indirect observation is a category of inverse problems that arises in numerous applications, such as medical imaging, geophysical prospection and nondestructive testing.

In this paper, we focus on transient models and the appearance of unknown nonlinearities. In order to find a finite dimensional representation of the latter, we make use of the powerful approximation properties and computational efficiency of \emph{neural networks (NNs)}. Due to the inherent ill-posedness of these inverse problems, however, regularization must be employed. We therefore study two such regularization methods: the variational Tikhonov method and the iterative projected Landweber method. These reconstruction methods are analyzed in the spirit of regularization theory, with a discretization by neural networks, as well as other general discretization schemes. We must therefore investigate the interplay between noise level, regularization parameter and discretization error \blue{(approximation error in case of discretization by NNs). In the language of machine learning, this is the interplay between approximation error and optimization/estimation error, with the impact of ill-posedness and data noise additionally taken into consideration. The resulting convergence analysis hints at a dependence of the network size (discretization parameter) on the noise level. This constitutes one of the main contributions of the paper.}

We point to the fact that a regularization theoretical viewpoint for the training problems has already been taken in \cite{AspriKorolevScherzer20}, although there the focus is on linear problems. \blue{In \cite{AspriKorolevScherzer20}, the authors solve $Au=y$ without making use of the linear forward map $A$, relying solely on the input-output training pairs $[u_m,y_m]_{m=1\ldots M}$ satisfying $Au_m=y_m$. This data-driven approach is interpreted as regularization by projection, where the subspaces are spanned by the training data. Along this line, \cite{BurgerEngl20} investigates the supervised training problem of  approximating a smooth function via one-layer feed-forward networks with noisy data as an ill-posed problem. This is shown to be equivalent to least-squares collocation for a linear integral equation. The core result is the derivation of an optimal choice of the network size depending upon on the data error $\delta$. In the same spirit, our work focuses on the connection between machine learning and regularization. Our objective is to establish a convergence analysis of regularization methods under the influence of the network approximation error in a nonlinear PDE. As the considered problem already exhibits multi-faceted complexity, namely parameter identification, a nonlinear model and unsupervised training, the task of deriving convergence rates as in \cite{BurgerEngl20} is deferred to future research.}

Our use of neural networks in parameter identification is inspired by \cite{DongHintermuellerPapafitsoros20}. There, the focus is on stationary problems, \blue{and the nonlinearity is represented by a neural network. In \cite{DongHintermuellerPapafitsoros20}, the network is learned beforehand via supervised training, and then it is inserted into the PDE model underlying the parameter identification. The supervised learning thus requires exact and full measurements of the state $u$, as well as physical parameters of the PDE, in order to form the training pairs. In contrast to this, we here consider time-dependent models; in addition, by virtue of our \emph{all-at-once formulation}, the supervised training is skipped, hence no access to the exact state and physical parameters is involved. The application of an all-at-once approach for  unsupervised learning is another main contribution of our paper.}

\blue{Another advantage of the all-at-once formulation lies in the fact that it avoids the evaluation of a parameter-to-state map, thus bypassing the need for a nonlinear PDE solver in a practical implementation. Additionally, this setting simplifies verification of the so-called tangential cone condition, a requirement for convergence guarantees of gradient-based methods, whose verification in many applications is neglected. This, may be considered another advantage of our approach.}

\blue{This work is a continuation of \cite{Nguyen:21}, in which we also parametrize the unknown nonlinearities in time-dependent PDEs by NNs. There, a so-called \emph{learning-informed parameter identification} was investigated by way of discretized inverse problems (i.e. when $f$ is already approximated by some NN). The present study develops a theoretical framework for \cite{Nguyen:21}, in the sense that we show convergence of the regularized and discretized reconstructions towards a ground truth. The approximation/discretization is incorporated into the regularization, which places a stronger emphasis on regularization theory in the light of existing literature. The analysis applies not only to discretization by NNs, but also to more general discretization schemes.}

 \blue{The field of deep learning for PDEs is well developed, with many novel results and techniques available in the literature. One such technique is \emph{physics informed neural networks} (PINNs) \cite{Raissi}, where one parametrizes the solution to the PDE, as opposed to using NNs to parametrize the unknown nonlinearities as in our case. A theoretical justification for using NNs to parametrize PDE solutions or parameter-to-solution maps can be found in \cite{Kutyniok22}.}
 
 \blue{Recovery of hidden physics laws from empirical observations is, in fact, an active field with a significant history. Recently, the rapid advances in computing power and data acquisition open the door for advanced techniques. For example \cite{Brunton16,Schaeffer} are concerned with the recovery of the governing PDE from full measurements of the state $u$. These two papers suggest to first construct a rich library of possible basis elements and then optimize the corresponding coefficients using sparse regression. Adding deep learning techniques, \cite{Raissi18} proposes to use two deep neural networks, one representing the solution $u$, the other  representing the nonlinear dynamics $f:=u_t-\Nc(t,x,u,\nabla u,\nabla u^2\ldots)$. Algorithmic differentiation is employed for computing the required derivatives. On the other hand, PDE-NET \cite{Long18} represents another flexible framework, in which one approximates the model $f$ by a feed-forward NN, while numerically approximating the differential operators $\nabla, \nabla^2\ldots$ by convolutional NNs. In all these mentioned studies, the problems are studied in a discrete setting and the collocation points $(t_i,x_i)$ range over the entire time and space $[0,T]\times\Omega$. In this work, we put more  emphasis on the aspect of the PDE model being derived from physical laws, and use data-driven methods solely to complement this approach. Further, our study combines a functional setting for the unknown physical parameters and states with a network parametrization for the unknown nonlinearity in an hybrid form.}

While \cite{AspriKorolevScherzer20,BurgerEngl20,DongHintermuellerPapafitsoros20, Raissi,Brunton16,Schaeffer,Raissi18,Long18,Nguyen:21} are the recent publications that our work is most closely related to, there is clearly a vast amount of existing and emerging literature on the mathematics of machine learning. In the context of data-driven inverse problems, we refer to \cite{ArridgeMaassOektemSchoenlieb:19} for an excellent review. \blue{For a profound exposition on the theory of deep learning, we refer to the lecture series and associated upcoming publication \cite{BookKutyniok}.}

\subsection{The inverse problem}

Quite often, the nonlinearity is not the only unknown quantity, but rather must be determined alongside other coefficients in the PDE, as exemplified in the following application.
\paragraph{Application.} 
Consider the problem of recovering the unknown nonlinearity $f$, 
the potential $c$, the source \blue{$\varphi$}, and the initial data $u_0$
in
\begin{equation}\label{IP_c-problem}
\begin{aligned}
&\dot{u} -\Delta u + c u +h(u)= f(u)+\source \quad&&\mbox{ in } \blue{\Omega\times}(0,T)\\
&u(0)=u_0 &&\mbox{ on } \blue{\Omega}\,,
\end{aligned}
\end{equation}
from measurements $y$ of the state $u$. \blue{The state $u$ is a function on the finite time line $(0,T)$ and the bounded, smooth domain $\Omega$ with its time derivative being denoted by $\dot{u}$}. In \eqref{IP_c-problem}, $h$ is the known nonlinear part of the model; the unknown nonlinear part $f$, which needs to be determined, plays the role of a model correction, thus helping to refine the physical model. The additional data $y$ available to identify the unknown quantities are observations of the state $u$ expressed via some observation operator $M$ (which could e.g., be the trace of $u$ at the boundary over time or its values at some fixed times instance(s) in $\Omega$)
\begin{equation}\label{Muy}
Mu=y\,.
\end{equation}
In more complex settings, for the  purpose of identifying the unknown functions, several repeated or possibly also different observations $y^{\indmeas}=M^{\indmeas}u^{\indmeas}$, $\indmeas=1\ldots K$ will be needed. These observations entail a variation in the data, and possibly also in some of the unknown coefficients, while the nonlinearity remains the same. Different measured data $y^{\indmeas}$ correspond to the model \eqref{IP_c-problem} at different parameters, thus at different states, i.e. $(c,\varphi,u_0)^{\indmeas}, u^{\indmeas}, \indmeas=1\ldots K$ may vary between observations, while the unknown function $f$ describing the underlying physical law is fixed. \blue{There are several real life inverse problems obeying this setting. In medical imaging MRI, this situations appears when different patients are scanned, resulting in $K$ sets of patient-dependent physical/body parameters. These patient-specific parameters, however, enter the same model governed by the same underlying physical law, e.g. the Bloch-Torrey equation model \cite{BenningEhrhardt}. Thus, the unknown nonlinear response $f$ can be considered as being fixed, while relaxation and diffusion parameters are allowed to vary between patients.}

\begin{remark}[\blue{uniqueness}]
\blue{In the context of restricted measurements \eqref{Muy} such as boundary observations, as relevant in tomographic applications, the question arises whether $f$ can be determined uniquely from these observations. Answers to this question can be found in the literature in the case of unknown $f$ and known initial data $u_0$ as well as coefficients $c$, $\varphi$, see, e.g., \cite{CannonDuChateau:1998,DuChateauRundell:1985,PilantRundell:1986}, or in the case of known $f$ and unknown initial data $u_0$ or coefficients $c$, $\varphi$, see, e.g., \cite{Isakov:1990,Isakov:2006}, but more rarely on simultaneous identifiability of all these quantities.
As the unknown $f$ case is the most relevant for our study, we point to the fact that a range condition on the exact state $u_{\textup{exact}}$
\[
u_{\textup{exact}}(\Omega\times(0,T))\subseteq
u_{\textup{exact}}(\omega\times(0,T)) 
\]
is essential for establishing unique recovery of $f(u)$ from observations of $u$ in some subset $\omega$ of $\Omega$ or its boundary. }
\end{remark}

\medskip

We will study such inverse problems in a more general framework of the following form.
\paragraph{Inverse Problem.} 
Before stating the inverse problem, we point out that in the model \eqref{IP_gen} below, the unknown nonlinearity $f:\R^{n+1}\to\R$ is identified with the corresponding \emph{Nemitskii operator} (cf. Section~\ref{sec:prelim})
\begin{align}\label{Nem-f}
f: \mathbb{R}^{n}\times \mathcal{V}\to \mathcal{W} \quad\text{via}\quad [f(\alpha,u)](x,t):=f(\alpha,u(x,t)),
\end{align}
and similarly, the known part $F:(0,T)\times X\times V\to W$ is identified with the Nemitskii operator 
\begin{align}\label{Nem-F}
 F:X\times \Vc\to \Wc \quad\text{via}\quad [F(\lambda,u)](t):=F(t,\lambda,u(t)),
\end{align}
with $\Vc$ denoting the state space and $\Wc$ denoting the image space of the model. Here, $\mathcal{V}\subseteq L^2(0,T;V)$ and $\mathcal{W}\subseteq L^2(0,T;W)$ are \emph{Bochner spaces}  (cf. Section~\ref{sec:prelim}) and with $V$ being a space of $x$ dependent functions, in \eqref{Nem-f} we make the identification $u(x,t)=(u(t))(x)$.

We now investigate the inverse problem of determining the physical parameters $\lambda^{\indmeas}\in X$ (parameter space), $u_0^{\indmeas}\in H$ (initial data space), $\alpha^{\indmeas}\in \mathbb{R}^n$, $\indmeas\in\{1,\ldots,K\}$, and the nonlinearity $f\in \mathcal{C}\subseteq C(\mathbb{R}^{n+1};\mathbb{R})$, in the evolution system
\begin{equation}\label{IP_gen}
\begin{aligned}
&\dot{u}=F(\lambda,u)+f(\alpha,u)\quad\mbox{ in }(0,T)\\
&u(0)=u_0
\end{aligned}
\end{equation}
from $K$ noisy measurements $y^{\indmeas,\delta}\in\mathcal{Y}^{\indmeas}$ (data space) of the states $u^{\indmeas}\in \mathcal{V}$ (state space) under the measurement operator $M^{\indmeas}$ according to
\begin{equation}\label{IP_gen_data}
\begin{split}
&y^{\indmeas} = M^{\indmeas} u^{\indmeas},  \quad M^{\indmeas}:\Vc\to\mathcal{Y}^{\indmeas}\\
&0\leq S^{\indmeas}(y^{\indmeas},y^{\indmeas,\delta})\leq\delta\quad \indmeas\in\{1,\ldots,K\},
\end{split}
\end{equation}
\blue{where $u^{\indmeas}$ solves \eqref{IP_gen} with $(\lambda,u_0,\alpha)=(\lambda^{\indmeas},u_0^{\indmeas},\alpha^{\indmeas})$.
}
Here, the distance between the exact data $y^{\indmeas}$ and the noisy data $y^{\indmeas,\delta}$ under the \blue{misfit} measure $S^{\indmeas}$ is assumed to be bounded by the noise level $\delta$. Typical choices of $S^{\indmeas}$ are norms (as in Section~\ref{sec:Landweber} below) or more general \blue{distance} measures such as the Kullback-Leibler divergence. If \eqref{IP_gen} represents an evolutionary PDE such as \eqref{IP_c-problem}, the spaces $V, W, X, H$ are typically Banach spaces of space-dependent functions, 
and $\Cc$ is a function space over $\mathbb{R}^{n+1}$ that will be approximated by neural networks later on.
The data spaces $\mathcal{Y}^{\indmeas}$ are Banach spaces as well and, depending on what type of observations are made, may consist of space and/or time dependent functions.

This parameter identification problem can equivalently be written as the \emph{all-at-once system} (cf. Section~\ref{sec:prelim})
\begin{equation}\label{aao}
\left(\begin{array}{l}E(\zeta^{\indmeas},u^{\indmeas},f)\\M^{\indmeas}u^{\indmeas}\end{array}\right)_{\indmeas=1}^K=\left(\begin{array}{l}0\\y^{\indmeas}\end{array}\right)_{\indmeas=1}^K, \qquad \zeta^{\indmeas}:=(\lambda^{\indmeas},u_0^{\indmeas},\alpha^{\indmeas})\in Z:=X\times H\times \mathbb{R}^n
\end{equation}
\blue{for the parameters $\zeta^{\indmeas}$, the nonlinearity $f$, and the states $u^{\indmeas}$,} with
\begin{equation}\label{aao-1}
\begin{split}
&E(\zeta,u,f):= \left(\begin{array}{l}\dot{u}-F(\lambda,u)\blue{-}f(\alpha,u)\\u(0)-u_0\end{array}\right)\in \mathcal{W}\times H.
\end{split}
\end{equation}

We denote by  $(\zeta^\dagger$, $u^\dagger$, $f^\dagger)$ an exact solution to the inverse problem, that is
\[
E(\zeta^\dagger,u^\dagger,f^\dagger)=0\,, \quad Mu^\dagger=y,
\]
\blue{leading to a vanishing PDE residual} and a perfect match of the measurements to the noise free data $y$.

In particular, the application \eqref{IP_c-problem} is a special case of \eqref{IP_gen} in the setting \eqref{aao} with
\begin{equation}\label{Flambda_c-problem}
\blue{K=1}, \quad n=0, \quad  \lambda=(c,\source), \quad
F(\lambda,u)=\Delta u - c u -h(u)+\source.
\end{equation}

\subsection{Preliminaries}\label{sec:prelim}

Before launching a detailed discussion of discretization for Tikhonov and Landweber regularization, we briefly elaborate on some concepts that have been mentioned in the preceding section.
\paragraph*{Bochner spaces.} 
Given a Banach space $V$, the \emph{Bochner space} $L^p(0,T;V)$ \cite[Section 1.5]{Roubicek} consists of the Bochner integrable functions $u:[0,T]\to V$ satisfying $\int_0^T\blue{\|u(t)\|^p_V} dt<+\infty$. It is a Banach space under the norm
\begin{align*}
&\|u\|_{L^p(0,T;V)}:=\left(\int_0^T\|u(t)\|^p_V\,dt\right)^{1/p} \qquad 1\leq p<\infty.
\end{align*}
Likewise, the Bochner spaces $L^\infty(0,T; V)$ and $C(0,T;V)$ are Banach spaces under the respective norms
\begin{align*}
&\|u\|_{L^\infty(0,T;V)}:=\sup_{t\in[0,T]}\|u(t)\|_V,
\qquad \|u\|_{C(0,T;V)}:=\max_{t\in[0,T]}\|u(t)\|_V.
\end{align*}
Given a convex Banach space $V_1$ and a locally convex Banach space $V_2\supset V_1$, we define the \emph{Sobolev-Bochner space} $W^{1,p,q}(0,T;V_1,V_2)$ \cite[Section 7.1]{Roubicek}, which itself is a Banach space, as
\begin{align*}
&W^{1,p,q}(0,T;V_1,V_2):=\{u\in \blue{L^p}(0,T;V_1):\dot{u}\in\blue{L^q}(0,T;V_2)\} \quad \quad 1\leq p,q\leq\infty,\\ &\|u\|_{W^{1,p,q}(0,T;V_1,V_2)}=\|u\|_{L^p(0,T;V_1)}+\|\dot{u}\|_{L^q(0,T;V_2)}.
\end{align*}
An  example that is used in Section \ref{Land-dis} is $W^{1,2,2}(0,T;V_1,V_2)=L^2(0,T;V_1)\cap H^1(0,T;V_2)$.

\paragraph*{Nemitskii operators.} A mapping  \blue{$\mathfrak{f}:I\times\spaceA\to\spaceB$ with Banach spaces $\spaceA$, $\spaceB$ and } $I\subset\R^d$ is called a \emph{Caratheodory mapping} if $\mathfrak{f}(\cdot,u)$ is measurable for all $u\in\blue{\spaceA}$ and $\mathfrak{f}(\blue{z},\cdot)$ is continuous for a.e. $\blue{z}\in I$. The so-called \emph{Nemitskii operator} $\mathfrak{F}$ assigns a function $v:I\to\blue{\spaceA}$ to a function $w:I\to\blue{\spaceB}$ by
\[[\mathfrak{F}(v)](\blue{z}):=\mathfrak{f}(\blue{z},v(\blue{z})).\]
In \eqref{Nem-f}, we have, \blue{for any fixed $\alpha\in\R^n$, $\mathfrak{f}=f(\alpha,\cdot)$, $v=u$,} $I=\Omega\times(0,T)$, \blue{$\spaceA=\spaceB=\R$}, while in \eqref{Nem-F}, we have \blue{for any fixed $\lambda\in X$, $\mathfrak{f}=F(\lambda,\cdot)$, $v=u$,} $I=(0,T)$, \blue{$\spaceA=V$, $\spaceB=W$}. 
\\
For a detailed discussion on Nemitskii operators in Bochner spaces, we refer to \cite[Sections 1.3, 1.4]{Roubicek}.

\paragraph*{All-at-once formulation.} The classical way to formulate the inverse problem \eqref{IP_gen}-\eqref{IP_gen_data} 
(for simplicity of exposition setting $K=1$ and therefore skipping the superscripts $\indmeas$)
is to construct the \emph{reduced} forward operator
\[G:Z\times\Cc\to\Yc\quad G(\zeta,f):=M\circ \partostate(\zeta,f)=y,\]
which composes the observation operator $M$ with the parameter-to-state map
\[\partostate:Z\times\Cc\to\Vc \quad \partostate(\zeta,f)=u,
\mbox{ where $u$ solves \eqref{IP_gen}}.\] 
This formulation involves evaluating well-definedness of $\partostate$ via unique existence theory for the nonlinear PDE \eqref{IP_gen}, and in practice requires solving this nonlinear equation. 

Alternatively, the \emph{all-at-once} approach formulates \eqref{IP_gen}-\eqref{IP_gen_data} into a system
\begin{align*}
E(\zeta,u,f)&=0\\
Mu&=y
\end{align*}
of model and observation equation as in \eqref{aao}-\eqref{aao-1}. Hence, we can define the forward operator
\[\mathbf{G}:Z\times \Vc\times \Cc\to\Wc\times\Yc \quad \mathbf{G}(\zeta,u,f):=(E(\zeta,u,f),Mu)=\blue{(0,y)}.\] 
The all-at-once formulation bypasses the construction of the parameter-to-state map $\partostate$, which is nonlinear and often requires restrictive assumptions on $F,f$. This formulation therefore allows more general classes of $F,f$, and is also advantageous in practical implementation, where a PDE solver is not needed. All-at-once approaches have been studied in PDE constrained optimization in \cite{KunischSachs,KupferSachs,LeibfritzSachs,LeHe16,orozco-ghattas-97b,shenoy-heinkenschloss-cliff-98,Taasan91} and more recently, for ill-posed inverse problems, in
\cite{BurgerMuehlhuberIP,BurgerMuehlhuberSINUM,HaAs01,aao16,KKV14b,LeHe16}; a comparison between the reduced and all-at-once formulation for time dependent problems can be found in \cite{Kaltenbacher:17,Nguyen:19}.

\paragraph*{Neural networks (NNs).}
In the setting of this paper, we make use of the \emph{feedforward neural network} of \emph{depth} $L$ on $(\alpha,u(x,t))\in\R^{n+1}$, expressed as a function \blue{of} the form
\[\Nc:\R^{n+1}\to\R \quad \Nc(\alpha,u(x,t)) :=A_L\circ\ldots A_1(\alpha,u(x,t)), \quad
A_\ell(\mathbf{z}) :=\sigma_\ell(w_\ell \mathbf{z}+b_\ell),\]
where the matrix $w_\ell\in \Lc(\R^{p_{\ell-1}}, \R^{p_\ell})$ and the vector $b_\ell\in \R^{p_\ell}$ are the so-called \emph{hyperparameters} at \emph{layer} $\ell=1 \text{ (input)}\ldots L \text{ (output)}.$ The \emph{activations} $\sigma_\ell:\R\to\R$ are nonlinear point-wise functions allowed to differ between layers, and $\sigma_L=\text{Id}$. In summary, at layer $\ell-1$ the affine operator $A_\ell$ transforms an input vector in $\R^{p_{\ell-1}}$ into one in $\R^{p_\ell}$, applies the activation $\sigma_\ell$ pointwise, and returns the input to the next layer $\ell$. Some standard activation functions include the RELU function $\sigma(z)=\max\{z,0\}$, tansig function $\sigma(z)=\tanh(z)$, softsign function $\sigma(z)=\frac{z}{1+|z|}$ and softplus function $\sigma(z)=\ln(1+e^x)$.

Based on the universal approximation theorem for smooth functions \cite{Hornik}, we use standard feedforward neural networks to approximate the nonlinearity $f:\R^{n+1}\to \R$ in the finite dimensional set $\Cc_N$ (cf. \eqref{CN}), whose number of hyperparameters is $N=\sum_{\ell=1}^L(p_{\ell-1}+1) p_\ell$ with $p_0=n+1, p_{L+1}=1$. Fitting this into the formulation of the inverse problem, we identify $f$ as a Nemitskii operator between Bochner spaces, as introduced in \eqref{Nem-f}.

\paragraph*{Notation}
\begin{itemize}[label=]

\item
We will use shortcut notations $L^2(L^2)=L^2(0,T;L^2(\Omega))$, $C(L^2)=C(0,T;L^2(\Omega))$, $L^2=L^2(\Omega)$, and analogously for some further Sobolev spaces, when they appear as subscripts in some norms or constants.
\item
We will make use of boundedness of some Sobolev embeddings according to, e.g., \cite[Chapter 4]{Adams}, \cite[Chapter 11]{Leoni:2009} and more generally denote embedding constants between spaces $X$ and $Y$ by $C_{X\to Y}$. Generic constants will be denoted by $C>0$, and continuity or compactness of embeddings is indicated by $X\embed Y$ or $X\compt Y$, respectively.
\item
Partial derivatives are denoted by subscripts, e.g., $f_\alpha$, $f_u$, while ordinary or total derivatives by a prime, e.g. $f'$.
\end{itemize}

The remainder of this paper \blue{is} organized as follows. In Sections~\ref{sec:Tikhonov}, we prove convergence of Tikhonov regularization with an appropriate choice of the regularization parameter. Section \ref{sec:Landweber} presents convergence results for Landweber regularization with an  appropriate stopping index. In both approaches, the discretization level $N$ needs to be chosen too, in order to achieve convergence as the noise level \blue{$\delta$} tends to zero. The required conditions are thoroughly discussed and interpreted for the particular Application \ref{IP_c-problem}.

\section{Tikhonov regularization}\label{sec:Tikhonov}

With positive definite model and data misfit as well as regularization functionals
\begin{align*}
&Q:\mathcal{W}\to[0,\infty] \qquad{s.t.}\quad Q(w)=0\Leftrightarrow w=0,\\
&S:Y^2\to[0,\infty] \qquad{s.t.}\quad S(y_1,y_2)=0\Leftrightarrow y_1=y_2,\\
&R_1:Z\times\mathcal{V}\to[0,\infty], \qquad R_2:\Cc\to[0,\infty],
\end{align*}
consider the objective functional $T_\gamma^\delta$ given by
\[
T_\gamma^\delta(\vec{\zeta},\vec{u}, f):=
\sum_{\indmeas=1}^K \Bigl(Q(E(\zeta^\indmeas,u^\indmeas,f))+ S(M^\indmeas u^\indmeas,y^{\indmeas,\delta}) + \blue{\gamma} R_1(\zeta^\indmeas, u^\indmeas)\Bigr) + \blue{\gamma} R_2(f).
\]
Here, $K$ is the number of parameters and states corresponding to $K$ different observations of the data $y^{\indmeas,\delta}$, while $f$ is the common nonlinearity across all experiments. \blue{The objective functional $T_\gamma^\delta$ depends on the noise level $\delta$, measured data $y^{m,\delta}$ and the regularization parameter $\gamma>0$.} We then define \emph{regularized approximations} as minimizers of $T_\gamma^\delta$, that is,
\begin{align}\label{Tikh_ori}
(\vec{\zeta}^{\gamma,\delta},\vec{u}^{\gamma,\delta}, f^{\gamma,\delta})
\in\mbox{argmin}_{(\zeta^1,u^1),\ldots,(\zeta^K,u^K)\in (Z\times\mathcal{V})^K, f\in \mathcal{C}} T_\gamma^\delta(\vec{\zeta},\vec{u}, f).
\end{align}

The unknown nonlinearity $f$ is approximated by NNs, that is, within
the finite dimensional set
\begin{align}\label{CN}
\mathcal{C}_N:=\{\mbox{neural networks on $\mathbb{R}^{n+1}$ with $N$ parameters}\}
\subseteq\mathcal{C}.
\end{align} 
Denoting by $N$ the discretization parameter, we define \emph{partially discretized regularized approximations} 
as
\begin{equation}\label{Tikh_N}
(\vec{\zeta}^{\gamma,\delta,N},\vec{u}^{\gamma,\delta,N}, f^{\gamma,\delta,N})
\in\mbox{argmin}_{(\zeta^1,u^1),\ldots,(\zeta^K,u^K)\in (Z\times\mathcal{V})^K, f\in \mathcal{C}_N} T_\gamma^\delta(\vec{\zeta},\vec{u}, f).
\end{equation}
\blue{In comparison to \eqref{Tikh_ori}, the discretization parameter $N$ in \eqref{Tikh_N} enters the minimization as another parameter, which needs to be properly controlled. The focus of this section is on deriving a rule for the regularization parameter $\gamma$ and the discretization parameter $N$ with respect to the noise level $\delta$, such that convergence of the Tikhonov regularization method is guaranteed.} 
For simplicity of exposition, we set $K=1$, and mention in passing that an alternative way to take into account multiple observations, as opposed to summing over them in the Tikhonov functional, is the use of Kaczmarz methods. That is, implementing a cyclic iteration over the individual observations, see, e.g., \cite{Nguyen:19} for the all-at-once setting relevant here, as well as the references therein. 

\bigskip

\subsection{Convergence}
We now study convergence of the Tikhonov regularized approximations in the sense of regularization, so as $\delta\to0$  
with an appropriate choice of regularizer parameter $\gamma(\delta)$ and discretization parameter $N(\delta)$.
\begin{assumption}\label{ass1}
There exist topologies $\blue{\tau_1:=}\tau_{Z\times\mathcal{V}}$ on $Z\times\mathcal{V}$ and $\blue{\tau_2:=}\tau_{\mathcal{C}}$ on $\mathcal{C}$ such that the following holds: 
\begin{enumerate}[label=(T\arabic*)]
\item\label{Tikh-ass-levelset} sublevel sets of $R_1$ are $\tau_{Z\times\mathcal{V}}$ compact, and sublevel sets of $R_2$ are $\tau_{\mathcal{C}}$ compact;
\item\label{Tikh-ass-closed} $(Q\circ E, S\circ M)$ is $\tau_{Z\times\mathcal{V}}\times\tau_{\mathcal{C}}$ sequentially closed:\\
$\forall (\zeta^j,u^j,f^j,y^j)_{j\in\mathbb{N}}\subseteq Z\times\mathcal{V}\times\mathcal{C}\times \mathcal{Y}\, :$
\[
\begin{aligned}
&\Bigl( (\zeta^j,u^j,f^j)\stackrel{\tau_{Z\times\mathcal{V}}\times\tau_{\mathcal{C}}}{\longrightarrow}(\bar{\zeta},\bar{u},\bar{f})
\mbox{ and }Q(E(\zeta^j,u^j,f^j))\to0
\\&\qquad 
\mbox{ and }S(M u^j,y^j)\to0
\mbox{ and }S(y,y^j)\to0 \Bigr)\\
&\Longrightarrow \Bigl(Q(E(\bar{\zeta},\bar{u},\bar{f}))=0\mbox{ and }S(M \bar{u},y)=0\Bigr);
\end{aligned}
\]
\item\label{Tikh-ass-lsc}  
$(Q\circ E, S\circ M)$, $R_1$, $R_2$ are $\tau_{Z\times\mathcal{V}}\times\tau_{\mathcal{C}}$ lower semicontinuous:
\\
$\forall (\zeta^j,u^j,f^j)_{j\in\mathbb{N}}\subseteq Z\times\mathcal{V}\times\mathcal{C}\, :$
\begin{align*}
&(\zeta^j,u^j,f^j)\stackrel{\tau_{Z\times\mathcal{V}}\times\tau_{\mathcal{C}}}{\longrightarrow}(\bar{\zeta},\bar{u},\bar{f}) \nonumber
\\
&\Longrightarrow \Bigl(
Q(E(\bar{\zeta},\bar{u},\bar{f}))\leq \liminf_{j\to\infty} Q(E(\zeta^j,u^j,f^j))
\mbox{ and }S(M\bar{u},y^\delta)\leq \liminf_{j\to\infty} S(Mu^j,y^\delta) \nonumber\\
&\qquad \quad\mbox{ and }R_1(\bar{\zeta},\bar{u})\leq \liminf_{j\to\infty} R_1(\zeta^j,u^j)
\mbox{ and }R_2(\bar{f})\leq\liminf_{j\to\infty} R_2(f^j)\Bigr).
\end{align*}
\end{enumerate}
\end{assumption}
\blue{The choice of  $R_1$ and $R_2$ is dictated by the continuity requirements \ref{Tikh-ass-closed} and \ref{Tikh-ass-lsc}. Accordingly, the topology $\tau_{Z\times\mathcal{V}}\times\tau_{\mathcal{C}}$ needs to be sufficiently strong. This topology is then linked to $R_1$, $R_2$ via the constraint on compactness of the sublevel sets expressed in \ref{Tikh-ass-levelset}. Overall, these assumptions are thus criteria to choose the regularizers $R_1$ and $R_2$.}
The latter two  requirements in \ref{Tikh-ass-lsc} are automatically satisfied if $R_1$, $R_2$ are defined by norms on $Z\times\mathcal{V}$ and $\mathcal{C}$, provided the spaces are reflexive or duals of separable spaces and $\tau_i$ is defined by the corresponding weak(*) topology.

\begin{proposition}\label{prop:Tikhonov_well-posed}
Under the assumptions \ref{Tikh-ass-levelset}-\ref{Tikh-ass-lsc}, the discrete minimization problems \eqref{Tikh_N} admit minimizers.
\end{proposition}
{\it Proof.}
The proof follows from standard results \cite{EHNBuch,Troeltzsch} that essentially assume compactness of \blue{sublevel sets} of $R_i, i=1,2$, $\tau_{Z\times\mathcal{V}}\times\tau_{\mathcal{C}}$ closedness of $(Q\circ E, S\circ M)$, and lower \blue{semicontinuity} of $T^\delta_\gamma$.\\
$\diamondsuit$

Thus, under these assumptions, the method is well defined by \eqref{Tikh_N}. In order to prove that it actually defines a convergent regularization method, we need further assumptions on the approximation quality and on the choice of the regularization and discretization parameters.

\blue{Moroever, we allow for \emph{inexact minimization} by introducing the tolerance $\eta\geq0$ in the relaxed definition 
\begin{equation}\label{Tikh_N_eta}
T_\gamma^\delta(\vec{\zeta}^{\gamma,\delta,N},\vec{u}^{\gamma,\delta,N}, f^{\gamma,\delta,N})
\leq T_\gamma^\delta(\vec{\zeta},\vec{u}, f) +\eta
\quad \forall 
(\zeta^1,u^1),\ldots,(\zeta^K,u^K)\in (Z\times\mathcal{V})^K, f\in \mathcal{C}_N.
\end{equation}
This definition actually does not even require existence of a minimizer.
}
 
\begin{assumption}\label{ass2}.
\begin{enumerate}[label=(T\arabic*)]
\setcounter{enumi}{3}
\item\label{Tikh-ass-approxNN} approximation by NNs:
\[
\begin{aligned}
&q_N:=\inf_{f_N\in \mathcal{C}_N} Q(E(\zeta^\dagger,u^\dagger,f_N))
+\gamma_N (R_2(f_N)-R_2(f^\dagger))\\
&\qquad=\inf_{f_N\in \mathcal{C}_N} (Q(E(\zeta^\dagger,u^\dagger,f_N))-Q(E(\zeta^\dagger,u^\dagger,f^\dagger)))
+\gamma_N(R_2(f_N)-R_2(f^\dagger))\to0\,,\\ 
&\mbox{ as }N\to\infty;
\end{aligned}
\]

\item\label{Tikh-ass-asymptotic} asymptotics of the parameters as $\delta\to0$: There exists $C>0$ such that
\[
\gamma(\delta)\to0\,, \quad \frac{\delta}{\gamma(\delta)}\leq C, \quad \frac{q_{N(\delta)}}{\gamma(\delta)}\leq C,
\blue{\quad  \frac{\eta}{\gamma(\delta)}\leq C,} 
\quad N(\delta)\to\infty, 
\blue{\quad  \eta(\delta) \to0} \
\mbox{ as }\delta\to0
\]
\end{enumerate}
\end{assumption}

\begin{proposition}\label{prop:Tikhonov_convergence}
Under Assumptions~\ref{ass1}, \ref{ass2} we have $\tau_{Z\times\mathcal{V}}\times\tau_{\mathcal{C}}$ subsequential convergence of \\
$(\zeta^{\gamma(\delta),\delta,N(\delta)},u^{\gamma(\delta),\delta,N(\delta)}, f^{\gamma(\delta),\delta,N(\delta)})$ to a solution of the inverse problem \eqref{aao}, \eqref{aao-1} as $\delta\to0$, \blue{i.e., every sequence $(\zeta^{\gamma(\delta_j),\delta,N(\delta_j)},u^{\gamma(\delta_j),\delta_j,N(\delta_j)}, f^{\gamma(\delta_j),\delta,N(\delta_j)})$ with $\delta_j\to0$ as $j\to\infty$ has a $\tau_{Z\times\mathcal{V}}\times\tau_{\mathcal{C}}$ convergent subsequence, and the limit of every $\tau_{Z\times\mathcal{V}}\times\tau_{\mathcal{C}}$ convergent subsequence solves the inverse problem.}

\end{proposition}

Note that unlike \cite{PoeResSch10} we assume convergence $q_N\to0$ of the ``discretization error'' to zero only at the exact solution, not uniformly over all elements of $Z\times\mathcal{V}\times\mathcal{C}$. Also, it is not necessary to assume any vector space structure on $\mathcal{C}_N$ and $S$ does not need to satisfy a triangle inequality.\\
{\it Proof.} By minimality, that is, 
$T_\gamma^\delta(\zeta^{\gamma,\delta,N},u^{\gamma,\delta,N}, f^{\gamma,\delta,N})\leq T_\gamma^\delta(\zeta,u,f_N)+\eta$ for all $\zeta\in Z$, $u\in\mathcal{V}$, $f_N\in\mathcal{C}_N$, setting $\zeta=\zeta^\dagger$, $u=u^\dagger$, \blue{and thus $Q(E(\zeta^\dagger,u^\dagger,f^\dagger)=0, Mu^\dagger=y$, shows that}
\blue{
\begin{equation}\label{estmin}
\begin{aligned}
T_\gamma^\delta&(\zeta^{\gamma,\delta,N},u^{\gamma,\delta,N}, f^{\gamma,\delta,N})\\
&\leq \inf_{f_N\in \mathcal{C}_N}\Bigl(Q(E(\zeta^\dagger,u^\dagger,f_N))+S(Mu^\dagger,y^\delta)+\gamma R_1(\zeta^\dagger,u^\dagger)+\gamma R_2(f_N)
+\eta\Bigr)\\
&
\leq Q(E(\zeta^\dagger,u^\dagger,f^\dagger))+S(y,y^\delta) +\gamma R_1(\zeta^\dagger,u^\dagger)+\gamma R_2(f^\dagger)
+\eta
\\&\qquad+\inf_{f_N\in \mathcal{C}_N}\Bigl(
\big(Q(E(\zeta^\dagger,u^\dagger,f_N)-Q(E(\zeta^\dagger,u^\dagger,f^\dagger))\big)
+\gamma \big(R_2(f_N)-R_2(f^\dagger)\big)\Bigr)\\
&=\delta+\gamma R_1(\zeta^\dagger,u^\dagger)+\gamma R_2(f^\dagger)
+q_N+\eta.
\end{aligned}
\end{equation}}
By employing $S\geq0$, $Q\geq0$, and dividing by $\gamma>0$, we then obtain
\begin{equation}\label{estR1R2}
R_1(\zeta^{\gamma,\delta,N},u^{\gamma,\delta,N})+R_2(f^{\gamma,\delta,N})\leq
R_1(\zeta^\dagger,u^\dagger)+R_2(f^\dagger)+\frac{q_N}{\gamma}+\frac{\delta}{\gamma}
\blue{+\frac{\eta}{\gamma}},
\end{equation}
which by \ref{Tikh-ass-levelset}--\ref{Tikh-ass-lsc} implies existence of a $\tau_{Z\times\mathcal{V}}\times\tau_{\mathcal{C}}$ convergent subsequence $(\zeta^j,u^j,f^j)_{j\in\mathbb{N}}$ of $(\zeta^{\gamma(\delta),\delta,N(\delta)},u^{\gamma(\delta),\delta,N(\delta)}, f^{\gamma(\delta),\delta,N(\delta)})_{\delta>0}$ with limit $(\zeta^*,u^*,f^*)$. 
Since from the same minimality estimate, by $R_1\geq0$, $R_2\geq0$ we also get 
\[
\begin{aligned}
&Q(E((\zeta^j,u^j,f^j)))+S(Mu^j,y^{\delta_j})\\
&\leq \gamma(\delta^j) \Bigl(R_1(\zeta^\dagger,u^\dagger)+R_2(f^\dagger)\Bigr)+q_N(\delta^j)+\delta^j
\blue{+\eta(\delta^j)} 
\, \to0\mbox{ as }j\to\infty
\end{aligned}
\]
for any such $\tau_{Z\times\mathcal{V}}\times\tau_{\mathcal{C}}$ convergent subsequence, from \ref{Tikh-ass-approxNN} we conclude that $(\zeta^*,u^*,f^*)$ solves the inverse problem $E(\zeta^*,u^*,f^*)=0$, $Mu^*=y$.
\begin{flushright}
$\diamondsuit$
\end{flushright}

Note that due to estimate \eqref{estR1R2} and $\tau_i$-lower semiconituity of $R_i$ in \ref{Tikh-ass-lsc},
if in addition to \ref{Tikh-ass-asymptotic} 
\[
\frac{\delta}{\gamma(\delta)}\to0, \quad \frac{q_{N(\delta)}}{\gamma(\delta)}\to0, 
\blue{\quad \frac{\eta(\delta)}{\gamma(\delta)}\to0, }
\quad \mbox{ as }\delta\to0,
\]
then the limit according to Proposition~\ref{prop:Tikhonov_convergence} is even an $R_1$, $R_2$ minimizing solution of the inverse problem, that is, $(\xi^*,u^*,f^*)=\min (R_1(\xi,u)+ R_2(f))$, where the minimum is taken over all $(\xi,u,f)$ solving the inverse problem. 

\begin{discussion}[Quantitative approximation error]\label{appox-rate}
\blue{Assumption \ref{Tikh-ass-approxNN} follows by application of the universal approximation theorem \cite{Hornik} in our setting. This theorem states that for any continuous function $f$ on a compact domain, there exists a NN with a sufficiently large number of neurons approximating $f$ with arbitrary prescribed accuracy. Recently, advanced studies on quantifying the size of NNs have been carried out, even in terms of width and depth, to obtain approximation rates. Seminal results in this direction include \cite{Barron93, Barron94}, which show an asymptotic approximation rate $\mathcal{O}(1/\sqrt{N}))$ in the $L^2$-norm of NNs with $N$ neurons and sigmoidal activation to any target function $f$ with finite Fourier moments. The study on approximation rates has greatly evolved in recent years \cite{Mhaskar93,Yarotsky17,Petersen18,CaoXieXu08}. For a full survey on approximation theory, we refer to  \cite{DeVore21}, as well as \cite[Section 1.4.2]{Kutyniok22}, \cite[Table 1]{Lu21} for brief summaries.}

\blue{Incorporating these approximation rates into $q_N$ in \ref{Tikh-ass-approxNN} enables an analysis for the convergence rate of Tikhonov regularization, under so-called source conditions, see, e.g., \cite{EHNBuch,SKHK:2012}. The asymptotics of the parameters \ref{Tikh-ass-asymptotic} shows that when $\delta\to0$, the NNs size should increase accordingly, that is, $N(\delta)\to\infty$. The relation reveals a choice of the network size dependent on the noise level $\delta$. This potentially reduces the overfitting problem caused by noisy training data. 
By virtue of \eqref{estmin}-\eqref{estR1R2}, the approximation errors both w.r.t the model $Q\circ E$ and w.r.t the regularizer $R_2$ contribute to the total convergence rate, hinting at a possible mutual effect of these two factors in the overall rate. 
If a convergence rate analysis for Tikhonov regularization integrating the quantitative approximation error can be carried out, a choice of $N$ with respect to $\delta$ (c.f. \cite{BurgerEngl20}), as well as $Q\circ E$ and $R_2$ can be made explicit.  We leave this interesting task for future research.}
\end{discussion}

\subsection{On Assumptions~\ref{ass1}, \ref{ass2}: Discussion, Examples, and Application}\label{sec:Tikh-dis}
\begin{discussion}[on \ref{Tikh-ass-levelset}-\ref{Tikh-ass-closed}]\label{Tikh-dis-closed}
\blue{Assume the following:}
\begin{enumerate}[label=(\roman*)]
\item \label{Tikh-dis-cont}
$\exists\, 
\tau_{\Wc\times H}
$
such that $E$ is $\tau_{Z\times\Vc}\times\tau_\Cc\too
\tau_{\Wc\times H}
$ continuous at the exact solutions of the PDE.
\item \label{Tikh-dis-S}
$\exists\, \tau_\Yc$ such that $M$ is $\tau_\Vc\too\tau_\Yc$ continuous at the exact states.
In addition, $S(y,y^j)\to 0$ implies $y^j\rightarrow y$ in $\tau_\Yc$.
\item \label{Tikh-dis-lsc}
$(Q,S)$ is $(
\tau_{\Wc\times H}
,(\tau_\Yc\times\tau_\Yc))$ lower semicontinuous.
\item 
$(R_1,R_2)$ is $(\tau_{Z\times\Vc},\tau_\Cc)$ lower semicontinuous and its sublevel sets are $(\tau_{Z\times\Vc},\tau_\Cc)$ compact.
\end{enumerate}
\blue{Then \ref{Tikh-ass-levelset}-\ref{Tikh-ass-closed} hold.}

First, $Q\circ E$ is lower semicontinuous as it  is a composition of a lower semicontinuous function and a continuous function assumed in \ref{Tikh-dis-cont}-\ref{Tikh-dis-lsc}.
Next, $Q\circ E$ is closed since by positivity of $Q$, lower semicontinuity of $Q\circ E$ and the premise of \ref{Tikh-ass-closed}, one has
\begin{align*}
&0\leq Q(E(\zb,\ub,\fb)) \leq \liminf_{j\to\infty}Q(E(\zeta^j,u^j,f^j)) \leq \lim_{j\to\infty}Q(E(\zeta^j,u^j,f^j)) = 0\\
&\text{thus}\quad Q(E(\zb,\ub,\fb)) =0.
\end{align*}
Note that closedness in the sense of \ref{Tikh-ass-closed} is weaker than in the standard definition, (see, e.g., \cite{PoeResSch10} and the references therein), as we require the closedness property only at the exact solutions of the PDE, i.e. at $Q(E(\zb,\ub,\fb)) =0.$\\
Furthermore, if $S(y,y^j)\to 0$, that is the premise of \ref{Tikh-ass-closed}, induces $y^j\rightarrow y$ in $\tau_\Yc$ (cf. \ref{Tikh-dis-S}), then $S\circ M$ is closed due to
\begin{align*}
0\leq S(M\ub,y) \leq \lim_{j\to\infty}S(Mu^j,y^j) =0 \quad 
\implies
\text{thus}\quad S(M\ub,y))=0,
\end{align*}
provided that $S$ is lower semicontinuous in its two arguments (see \ref{Tikh-dis-lsc}).
\end{discussion}

\begin{remark}[on \ref{Tikh-ass-lsc}]
In Discussion \ref{Tikh-dis-closed}, if $E$ is 
$\tau_{Z\times\Vc}\times\tau_\Cc\too\tau_{\Wc\times H}$
continuous on the whole space $Z\times\Vc$ and $M$ is $\tau_\Vc\too\tau_\Yc$ continuous on $\Vc$, then lower semicontinuity of $T_\gamma^\delta$  on $Z\times\mathcal{V}\times\mathcal{C}_N$ \eqref{Tikh-ass-lsc} holds. In some particular examples where convexity of $R_i, \blue{i=1,2}$ is given, e.g. $R_i=\|\cdot\|^p,\blue{ i=1,2}$, for some $p\in[1,\infty]$ weaker conditions on continuity of $E, M$ might be sufficient.
\end{remark}

\begin{remark}
In case of full measurement, the term $S(Mu,y)=S(u,y)$ can play the role of a regularizer on $u$ with $\tau_\Vc=\tau_\Yc$.
\end{remark}

\begin{discussion}[on \ref{Tikh-ass-approxNN}]
The topology $\tau_\Cc$ induced by $R_2$ could be chosen as the weak$^*$ topology induced by the $L^\infty$-norm to make use of available approximation rates of deep neural networks to smooth functions. In particular, these rates are with respect to arbitrary depths (number of layers) and widths (number of neurons per layer) \cite[Table 1]{Lu21} to which $N$ in \eqref{Tikh_N} generally refers.\\
This and the discretization error assumption \ref{Tikh-ass-approxNN}
\begin{align*}
q_N:=\inf_{f_N\in \mathcal{C}_N} Q(E(\zeta^\dagger,u^\dagger,f_N))
=\inf_{f_N\in \mathcal{C}_N} Q(E(\zeta^\dagger,u^\dagger,f_N(u^\dagger))) \to0
\end{align*}
require uniform boundedness only on the exact state $u^\dagger$. Therefore, a candidate for $R_2$ is $R_2(f_N)=\|f_N\|_\lomy$ with $\Omega_{y^\dagger}=u^\dagger((0,T)\times\Omega)$.
\end{discussion}

\medskip

In the following examples of settings satisfying Assumptions~\ref{ass1}, \ref{ass2}, we consider reflexive spaces or duals of separable spaces.
\begin{example}\label{Tikh-ex-norm}
Let
\begin{alignat*}{3}
&Q(E(\zeta,u,f_N))=\|E(\zeta,u,f_N)\|^2_{\Wc\times H}, \qquad && S(Mu,y)=\|Mu-y\|^2_\Yc,\\ &R_1(\zeta,u)=\|(\zeta,u)\|^2_{Z\times\Vc}, && R_2(f)=\|f\|^2_{W^{1,\infty}(\Omega_y)},
\end{alignat*}
with a bounded interval $\Omega_y\subset\mathbb{R}$ containing $\Omega_{y^\dagger}=u^\dagger((0,T)\times\Omega)$
as detailed below. Since the required compactness and continuity properties are straightforward on the finite dimensional space $\mathbb{R}^n$, for simplicity of exposition we skip $\alpha$ as an argument of $f$.

Then:
\begin{itemize}
\item 
Let $\tau_X\times\tau_\Vc$ be the weak topology on $X\times\Vc$ and assume that
\begin{equation}\label{condEtil-Tikh}
\begin{aligned}
\tilde{E}:=(\frac{d}{dt}-F)\text{ is }(X\text{ weak})\times(\Vc\text{ weak})\too (\Wc\text{ weak})\\\text{ continuous at the exact solution.}
\end{aligned}
\end{equation} 
This weak continuity, thus closedness \ref{Tikh-ass-closed}, depends on the PDE models and the choice of function spaces. 

\item
Now, with $\tau_\Cc$ being the weak topology on ${W^{1,\infty}(\Omega_y)}$, we show continuity of the rest of $E$. In particular, we prove that $(u_N,f_N)\to(u,f)$ in the topology $(\Vc \text{ weak})\times (W^{1,\infty}(\Omega_y) \text{ weak})$  implies $f_N(u_N)\to f(u)$ weakly in $\Wc$ under appropriate conditions on $\Vc$, $\Wc$ to be derived here. First, we observe
\begin{align}\label{NNs-derivative}
&f_N(u_N(x,t))-f(u(x,t))= \left( f_N(u_N(x,t))-f_N(u(x,t))\right)+\left(f_N(u(x,t))-f(u(x,t))\right)\nonumber\\
&= \int_0^1 (f_N)' \left(u(x,t)+\theta(u_N(x,t)-u(x,t))\right)d\theta\,(u_N(x,t)-u(x,t))\\
&\qquad\qquad +\left(f_N(u(x,t))-f(u(x,t))\right) \nonumber\\[1ex]
&\langle f_N(u_N(x,t))-f(u(x,t)),\psi\rangle_{\Wc,\Wc^*}\\
&\leq \|(f_N)'\|_{L^\infty(\Omega_y)}\|u_N-u\|_{L^p(\tom)}\|\psi\|_{L^{p^*}(\tom)}+\|f_N-\blue{f}\|_{L^\infty(\Omega_y)}\|\psi\|_{L^1(\tom)}, \nonumber
\end{align}
for any $\psi\in\Wc^*$,
with $p\in[1,\infty]$ and $p^*$ being the conjugate index of $p$.
If 
\begin{equation}\label{condV_Tikh}
\Vc\subset L^\infty(\tom),
\end{equation} 
then for $u_N\stackrel{\Vc}\rightharpoonup u$, one has $\|u_N\|_{L^\infty\tom},\|u\|_{L^\infty\tom}\leq C, \forall N\in\N$, and may set $\Omega_y:=[-C-1,C+1]$. Note that the inclusion $\Vc\subset L^\infty(\tom)$ allows us to apply the fact that neural networks are dense in the space of smooth functions on compact sets. Next, $f_N \rightharpoonup f$ in $W^{1,\infty}(\Omega_y)$ shows that $\|(f_N)'\|_{L^\infty(\Omega_y)}$ is bounded for all $N$, and due to $W^{1,\infty}(\Omega_y)\compt L^\infty(\Omega_y)$ we have $ f_N\to f$ in $L^\infty(\Omega_y)$. If 
\begin{equation}\label{condVW_Tikh}
\Vc\compt L^{p_W}(\tom)\subset\Wc
\end{equation}
for some $p_W\in[1,\infty]$, then \eqref{NNs-derivative} shows $f_N(u_N)\rightharpoonup f(u)$ in $\Wc$, meaning $\tau_\Vc\times\tau_\Cc\too\tau_\Wc$ continuity of $(u,f)\mapsto f(u)$ on $\Vc\times\Cc$.\\
Recall that for closedness of $(u,f)\mapsto f(u)$, we require only its continuity at exact solutions $(\ud,\fd)$ of the PDE. Therefore, by 
\[
\begin{aligned}
&f_N(u_N(x,t))-\fd(\ud(x,t))\\
&\qquad= \left(\fd(u_N(x,t))-\fd(\ud(x,t))\right)+\left( f_N(u_N(x,t))-\fd(u_N(x,t))\right)
\end{aligned}
\]
we only need to assume boundedness of $\|(\fd)'\|_{L^\infty(\Omega_y)}$ thus can choose the weaker $R_2(f)=\|f\|^2_\Cc,\Cc\compt L^\infty(\Omega_y)$. Note that the inclusions \eqref{condV_Tikh}, \eqref{condVW_Tikh} are still needed.

\item
The part on the initial condition $E_0(u_0,u)=u(0)-u_0$ is linear, thus requiring just the embedding $\Vc\hookrightarrow C(0,T;H)$ and the regularizer $R_1(u_0)=\|u_0\|_H$ induces $\tau_H$, the weak topology on $H$.

\item
Regarding the observation $M$, if $M$ is linear and bounded, then it is $\Vc$ weak $\too\Yc$ weak continuous.
\end{itemize}
\end{example}

\begin{example}[norm of the hyperparameter as $R_2$]\label{Tikh-ex-hyperpar}
In the previous example, we consider the Sobolev $W^{1,\infty}$-norm for $f,f_N\in\Cc$. As the discretized regularization is carried out for $f_N\in \Cc^N$, the space of neural networks of $N$ hyperparameters, a natural question is whether one can replace the Sobolev norm by some equivalent norm on the hyperparameters. The answer in the general case is no. Consider e.g. the function $f(x)=x$ expressed via a 2 layers neural network of identity activation $f(x)=x=1\cdot\text{Id}(x+b)-b$. The hyperparmeters are $\theta:=(w_2,b_2,w_1,b_1)=(1,b,1,-b)$. So, when $b$ tends to infinity $\|\theta\|\to\infty$ while $\|f\|_{L^\infty(\Omega)}<\infty$ for any bounded domain $\Omega$.\\
Let us study a standard case of a neural network with fixed depth two
\[f_N(y)=W^N_2\cdot\sigma(W^N_1y+b^N_1)+b^N_2,\] where $y\in\R, \sigma:\R\to\R, W^N_1\in\R^{N\times1}, b\in\R^{N}, W^N_2\in\R^{1\times N},b^N_2\in\R$, and $\cdot$ denotes matrix multiplication. This means when $N\to\infty$, the width of the neural network sequence tends to infinity. Assuming that $\sigma$ is Lipschitz continuous with Lipschitz constant $L_\sigma$, we have
\begin{align} \label{NNs-lip}
|f_N(u_N)-f_N(u)|\leq L_\sigma|W^N_2|\cdot|W^N_1||u_N-u|,
\end{align}
where $|\cdot|$ represents element-wise absolute value. The class of Lipschitz activations used in practice is large; some examples include ReLU $\sigma(x):=\max\{0,x\}$ (with approximation rates),  tansig $\sigma(x):=\tanh(x)$, softplus $\sigma(x):=\ln(1+e^x)$, sigmoid or soft step $\sigma(x):=\frac{1}{1+e^{-x}}$, softsign $\sigma:=\frac{x}{1+|x|}$ etc.\\
Furthermore, we assume that $\sigma$ is coercive in the sense that $\exists C_\sigma>0: |y|\leq C_\sigma|\sigma(y)|,\forall y\in\R$, $\sigma(y)\ge 0$ for $y\geq 0$, all hyperparameters are nonnegative. We can then estimate
\begin{align*}
&|W^N_2|\cdot|W^N_1|\leq |W^N_2|\cdot(|W^N_1|+b^N_1)\leq C_\sigma W^N_2\cdot\sigma(W^N_1+b^N_1) \\ 
&\leq C_\sigma \left(W^N_2\cdot\sigma(W^N_11+b^N_1)+b^N_2\right)
=C_\sigma f_N(1)\leq \sup_{x\in[0,1]} C_\sigma|f_N(x)|=C_\sigma\|f_N\|_{L^\infty([0,1])}.
\end{align*}
Combining this with \eqref{NNs-lip} and \eqref{NNs-derivative}, we can replace boundedness of $\|(f_N)'\|_{L^\infty(\Omega)}$ by boundedness of $\|f_N\|_{L^\infty(\Omega_y)}$, and use the weaker regularizer $R_2=\|\cdot\|^2_{L^\infty(\Omega_y)}$ instead of $\|\cdot\|^2_{W^{1,\infty}(\blue{\Omega_y})}$. When considering $\R^+$, some examples for Lipschitz continuous and coercive activation functions are: ReLU, Leaky ReLU (coercive on $\R$), softplus etc. 
Assume further that the exact $\fd$ can be expressed exactly via a neural network, possibly with infinitely many hyperparameters, say $\fd\in\Cc_\infty$ with $\fd(0)=0$, similar to \eqref{NNs-lip} we have
\begin{align*}
&|\fd(y)|=|f^\dagger(y)-\fd(0)|\leq L_\sigma|W_2|\cdot|W_1||y|\leq  |\Omega_y|L_\sigma |W_2|\cdot|W_1|\leq C \|W_2\|_{\ell^2} \|W_1\|_{\ell^2} \leq C\|\theta\|_{\ell^2}^2.
\end{align*} 
Then one can also use the stronger norm $\|\cdot\|_{\ell^2}$ in the regularizer $R_2$, \blue{ alternatively $\|\cdot\|_{\ell^1}$, due to norm equivalence in finite dimensional hyperparameter spaces. The application of sparsity-promoting techniques, such as incorporating $\ell^1$ regularizers, has been proven as one of the remedies for overfitting in machine learning in practice. Indeed, the sparse optimization performs feature selection, yielding more interpretable trained models \cite{Tibshirani}.}
\end{example}

\begin{example}\label{Tikh-ex-TV} Consider $R_1(u)=\text{TV}(u),$ the total variation of $u$ on $\Omega_T:=(0,T)\times\Omega, \Omega\subset\R^{d}$.

In order for $R_1$ to be lower semicontinuous and have $\tau_\Vc$ compact sublevel sets, we have some options:
\begin{enumerate}
\item 
$\tau_\Vc$ is the weak$^*$ topylogy on BV$(\Omega_T)$, the space of functions of bounded variation on $(0,T)\times\Omega$. Recall that $u^j\stackrel{*}{\rightharpoonup}u$ in BV$(\Omega_T)$ is defined as $u^j\stackrel{L^1}{\to}u, \text{TV}(u^j)\to \text{TV}(u)$. TV is weak$^*$ lower semicontinuous on BV$(\Omega_T)$, and its sublevel sets are weak$^*$ compact in BV$(\Omega_T)$ \cite{BurgerOsher2013}.
\item
$\tau_\Vc$ is the strong topology on $L^{\frac{d+1}{d}-\epsilon}(\Omega_T)$ for arbitrary small $\epsilon>0$. Application of the compact embedding $\text{BV}(\Omega_T)\compt L^{\frac{d+1}{d}-\epsilon}(\Omega_T)$ yields that TV is lower semiconinuous on $L^{\frac{d+1}{d}-\epsilon}(\Omega_T)$, and its sublevel sets are compact in $L^{\frac{d+1}{d}-\epsilon}(\Omega_T)$ \cite[Theorem 2.5]{AcarVogel94}.
\item
$\tau_\Vc$ is the weak topology on $L^\frac{d+1}{d}(\Omega_T)$. Weak compactness of the sublevel sets is clear from the compact embedding mentioned above. Weak lower semicontinuity of TV was shown, e.g. in \cite[Theorem 2.3]{AcarVogel94}.
\end{enumerate}

Let us consider, for instance, the second case where $\tau_\Vc$ is the strong topology on $\Vc=L^{\frac{d+1}{d}-\epsilon}(\Omega_T)$ with $\tilde{E}=\dot{u}-\Delta u$, and assume $u^j\stackrel{\Vc}\to u, u^j(T)\stackrel{L^1(\Omega)}\to u(T),u^j(0)=u(0)=0, u^j(\partial\Omega)=u(\partial\Omega)=0$. Let $\tilde{\epsilon}=\epsilon d^2/(1-\epsilon d)$, then due to the estimate
\begin{align*}
&\langle \dot{u}-\dot{u}^j-\Delta (u-u^j),\psi\rangle_{\Wc,\Wc^*}=\int_{\Omega_T}(u-u^j)(-\dot{\psi}-\Delta\psi)\,dx\,dt + \int_\Omega (u-u^j)(T)\psi(T)\,dx\\
&\leq\|u-u^j\|_{L^{\frac{d+1}{d}-\epsilon}(\Omega_T)} \|\dot{\psi}+\Delta\psi\|_{L^{d+1+\tilde{\epsilon}}(\Omega_T)} \\
&\qquad+ C_{W^{1,d+1+\tilde{\epsilon}}(\Omega)\to L^\infty(\Omega)}\|u(T)-u^j(T)\|_{L^1(\Omega)}\|\psi(T)\|_{W^{1,d+1+\tilde{\epsilon}}(\Omega)},
\end{align*}
for any $\psi\in\Wc^*$,
one can chose $\tau_\Wc$ as the strong topology on 
\[
\Wc:=\left(L^{d+1+\tilde{\epsilon}}(0,T;W^{2,d+1+\tilde{\epsilon}})\cap W^{1,d+1+\tilde{\epsilon}}(0,T;L^{d+1+\tilde{\epsilon}}(\Omega))\right)^*.
\] 
Note that \blue{continuity of the embedding} $\Wc^*\embed C(0,T;W^{1,d+1+\tilde{\epsilon}}(\Omega))$ \cite[Lemma 7.3]{Roubicek} \blue{implies finiteness of} $\|\psi(T)\|_{W^{1,d+1+\tilde{\epsilon}}(\Omega)}$.

Regarding $f$, since $\Vc \not\subset L^\infty(\Omega_T)$, in order to obtain uniform boundedness of $u^j,u$, we invoke full measurement data in a sufficiently strong observation space, e.g. $M=\text{Id}, \Yc=L^\infty(\Omega_T)$. Then observe that the inclusions $\Vc\compt L^{\frac{d+1}{d}-\epsilon}(\Omega_T)\subset \Wc$ hold, so convergence of the  neural network sequence $f^j(u^j)\stackrel{j\to\infty}\to \fd(\ud)$, as discussed in Example \ref{Tikh-ex-norm}, is guaranteed.

As such, we have two types of convergence for the sequence $u^j$: the strong convergence in $\Vc=L^{\frac{d+1}{d}-\epsilon}(\Omega_T)$, and the weak$^*$ convergence in $\Yc=L^\infty(\Omega_T)$. These types of convergence are in general not equivalent. An example for this is the sequence of Rademacher functions $f_n:[0,1]\to\{0,1\}$ \cite[Example 4.13]{PoeResSch10}
\[f_n(x)=(-1)^{i+1} \quad \text{for}\quad x\in[(i-1)/2^n, i/2^n], 1\leq i\leq 2^n,\]
which weak$^*$ converges to zero in $L^\infty([0,1])$, but not in the $L^1$-norm, thus not in the $L^{\frac{d+1}{d}-\epsilon}$-norm.
\end{example}

\begin{example}
Consider $S=\text{KL}$, the Kullback–Leibler divergence  defined by 
\begin{align*}
\text{for } y\in L^1(\Omega_T), \text{KL}(g,y):=
\begin{cases}
\int_{\Omega_T} y\left(\dfrac{g}{y}-\log\left(\dfrac{g}{y}\right)-1\right)\,dx\,dt \qquad & g,y\geq 0\text{ a.e}.\\
\infty &\text{else}
\end{cases}
\end{align*}
It is clear that $S=\text{KL}$ does not satisfy a triangle inequality, a situation that is taken into account in this work. Positivity of KL is obvious as $(g/y-1)\geq\log(g/y)$ and, KL$(g,y)=0$ iff $g=y=0$.

\cite[Lemma A.2]{BrediesHoller20} states that $\text{KL}(y,y^j)\to 0$ implies $\|y^j-y\|_{L^1(\Omega')}\to 0$ 
for some $\Omega'$, and for $\{g^j\}\in L^1(\Omega')$ with $g^j\rightharpoonup g$ in $L^1(\Omega')$ as $j\to\infty$, then $\text{KL}(g,y^j)\leq \liminf_{j\to\infty}\text{KL}(g^j,y^j)$. Fitting into our framework, in particular for existence of $\tau_\Yc$  in Discussion \ref{Tikh-dis-closed}, from $S(y,y^j)=\text{KL}(y,y^j)\to 0$ inducing $y^j\to y$ in $L^1(\Omega_T)$, one can choose $\tau_\Yc$ as the strong topology on $L^1(\Omega_T)$. Also by this lemma, $S$ is $\tau_\Yc\times\tau_\Yc$ lower semicontinuous. Therefore, we need $M$ to be $\tau_\Vc\too L^1(\Omega_T)$ continuous; this condition is very much obtainable in practice. In case $M=\text{Id}$, an estimate similar to the one in  Example \ref{Tikh-ex-TV} could be carried out for $u^j\to u$ in $L^1(\Omega_T)$. Still, convergence of the \blue{neural} network part requires $R_1(u)=\|u\|_\Vc$ with $\Vc\subset L^\infty(\tom)$ and $L^1(\tom)\subset\Wc$.
\end{example}

\paragraph*{Application.} We now return to Application \eqref{IP_c-problem}, \eqref{Muy} and from Propositions~\ref{prop:Tikhonov_well-posed} and \ref{prop:Tikhonov_convergence} conclude a result for Tikhonov regularization in the setting of Example~\ref{Tikh-ex-norm}
\begin{equation}\label{Tikh_Appl}
(c^{\gamma,\delta,N},\source^{\gamma,\delta,N},u_0^{\gamma,\delta,N},u^{\gamma,\delta,N}, f^{\gamma,\delta,N})
\in\mbox{argmin}_{(c,\source,u_0,u,f_N)\in X_c\times X_\source\times H\times\Vc\times\mathcal{C}_N} T_\gamma^\delta(c,\source,u_0,u,f_N)
\end{equation}
for 
\begin{equation}\label{Tikh_Appl_func}
\begin{aligned}
T_\gamma^\delta(c,\source,u_0,u,f_N)=&
\|\dot{u}-\Delta u + cu +h(u)-\source-f_N(u)\|^2_{\Wc\times H}+ \|Mu-y\|^2_\Yc \\
&+ \blue{\gamma} \|(c,\source,u_0,u)\|^2_{X_c\times X_\source\times H\times\Vc} + 
\blue{\gamma} \|f_N\|^2_{W^{1,\infty}(\Omega_y)}\,,
\end{aligned}
\end{equation}
where we can replace the ${W^{1,\infty}(\Omega_y)}$ norm of $f_N$ by the hyperparameter norm according to Example~\ref{Tikh-ex-hyperpar}.

For this purpose, recall the following requirements on the underlying spaces: \eqref{condEtil-Tikh}, \eqref{condV_Tikh}, \eqref{condVW_Tikh}, as well as boundedness of $M:\Vc\to\Yc$ and of $\mbox{tr}_{t=0}:\Vc\to H$. Here,  we have the operator
$\tilde{E}u=\dot{u}-\Delta u + cu +h(u)-\source$, and the $X$ space is decomposed as $X=X_c\times X_\source$; recall that $h$ is known.
We use the spaces
\begin{equation}\label{W_Appl-Tikh}
\begin{aligned}
&H=W^{t_V,q_V}(\Omega)\,, \quad X_c=L^r(\Omega), \quad 
X_\source=\Wc=W^{-s,p}(0,T;W^{-t,q}(\Omega))\,, \\
&\Vc=W^{1-s,p}(0,T;W^{-t,q}(\Omega))\cap W^{-s,p}(0,T;W^{2-t,q}(\Omega))\cap W^{s_V,p_V}(0,T;W^{t_V,q_V}(\Omega)) 
\end{aligned}
\end{equation}
with 
\begin{equation}\label{sptqsVpV}
s_V>\frac{1}{p_V}\,, \quad t_V>\frac{d}{q_V}\,, \quad r\leq p_W:=\min\{p,q\}
\end{equation}
to satisfy \eqref{condV_Tikh}, \eqref{condVW_Tikh} and part of \eqref{condEtil-Tikh}. To see the latter for the $c$ part of the operator $\tilde{E}$, observe that for any sequence $(c_n,u_n)$ converging weakly to $(c,u)$ in $X_c\times\Vc$, by our choice of $s_V,p_V,t_V,q_V$ there exists a subsequence $(c_{n_k},u_{n_k})$ such that $c_{n_k}$ converges weakly in $L^r(\Omega)$ and $u_{n_k}$ converges strongly in $L^\infty(0,T;L^\infty(\Omega))$, so that for any $\psi\in\Wc^*\subseteq L^{p_W^*}(\tom)$ we have $\int_0^T\psi u\in L^{r^*}(\Omega)$ and thus 
\[
\int_0^T\int_\Omega \bigl(c_{n_k} u_{n_k} -cu\bigr)\psi \, dx\, dt=
\int_0^T\int_\Omega c_{n_k} (u_{n_k}-u)\psi\, dx\, dt 
+ \int_\Omega (c_{n_k}-c)\int_0^Tu\psi\, dt\, dx\to0\,. 
\]
Likewise, it is straightforward to see that on the strength of the embeddings available for $\Vc$, it suffices to assume continuity of the real function $h$ to achieve 
$(\Vc\text{ weak})\too (\Wc\text{ weak})$ continuity of the mapping $u\mapsto h(u)$ contained in $\tilde{E}$. 
Note that continuity of $\mbox{tr}_{t=0}:\Vc\to H$ also holds for any subspace $H$ in which 
$W^{t_V,q_V}(\Omega)$ is continuously embedded.
\begin{corollary}
With the spaces according to \eqref{W_Appl-Tikh}, \eqref{sptqsVpV}, $M:\Vc\to\Yc$ linear and bounded as well as $h:\mathbb{R}\to\mathbb{R}$ continuous, Tikhonov regularization is well-defined by \eqref{Tikh_Appl}, \eqref{Tikh_Appl_func}, and with a choice of $\gamma(\delta)$, $N(\delta)$ according to Assumption~\ref{ass2}~\ref{Tikh-ass-asymptotic} we have weak(*) subsequential convergence  of 
$(c^{\gamma(\delta),\delta,N(\delta)},\source^{\gamma(\delta),\delta,N(\delta)},u_0^{\gamma(\delta),\delta,N(\delta)},u^{\gamma(\delta),\delta,N(\delta)}, f^{\gamma(\delta),\delta,N(\delta)})$ to a solution of the inverse problem 
\eqref{IP_c-problem}, \eqref{Muy} as $\delta\to0$.
\end{corollary}

\def\bfw{\mathbf{w}}
\def\bfx{\mathbf{x}}
\def\bfy{\mathbf{y}}
\def\bfz{\mathbf{z}}
\def\bfF{\mathbb{F}}
\def\bfM{\mathbb{M}}
\def\bfX{\mathbb{X}}
\def\bfY{\mathbb{Y}}
\def\Hadj{\ast}
\def\altil{\tilde{\alpha}}
\def\betil{\tilde{\beta}}
\def\latil{\tilde{\lambda}}
\def\mutil{\tilde{\mu}}
\def\ftil{\tilde{f}}
\def\gtil{\tilde{g}}
\def\htil{\tilde{h}}
\def\util{\tilde{u}}
\def\vtil{\tilde{v}}
\def\wtil{\tilde{w}}
\def\ytil{\tilde{y}}
\def\ztil{\tilde{z}}
\section{Landweber iteration}\label{sec:Landweber}
In this section, for simplicity of exposition, 
collecting all unknowns 
$\lambda^{\indmeas},u_0^{\indmeas},\alpha^{\indmeas},u^{\indmeas},f$ in a single variable $\bfx$
and setting $\bfy=\left(\begin{array}{l}0\\y^{\indmeas}\end{array}\right)_{\indmeas=1}^K$
we rewrite \eqref{aao} as an operator equation
\begin{equation}\label{Fxy}
\bfF(\bfx)=\bfy.
\end{equation}
Moreover, we restrict the setting to Hilbert spaces $\bfX$, $\bfY$ with Hilbert space adjoints denoted by a superscript $\ast$.

Landweber iteration defines regularized approximations as gradient descent steps for the least squares cost functional $\|\bfF(\bfx)-\bfy^\delta\|^2$, explicitly,
\[
\bfx_{k+1}^\delta=\bfx_k^\delta-\bfF'(\bfx_k^\delta)^\Hadj(\bfF(\bfx_k^\delta)-\bfy^\delta)\,.
\]
Here, the stopping index $k=k_*(\delta,\bfy^\delta)$, which depends on the noise level $\delta$ and data $\bfy^\delta$, acts as a regularization parameter. In order to accommodate additional constraints, e.g., on the magnitude or sign of $\bfx$, we consider a subset $\bfM$ of $\bfX$. Constraints formulated by membership in the subset $\bfM$ can be incorporated by projection via
\begin{equation}\label{LWproj}
\bfx_{k+1}^\delta= P_\bfM \Bigl(\bfx_k^\delta-\bfF'(\bfx_k^\delta)^\Hadj(\bfF(\bfx_k^\delta)-\bfy^\delta)\Bigr),
\end{equation}
where the metric projection operator $P_\bfM$ onto a closed convex set $\bfM$ is characterized by the variational inequality
\begin{equation}\label{proj}
\bfx=P_\bfM(\tilde{\bfx})\ \Leftrightarrow \ \left(\bfx\in \bfM \mbox{ and }
\forall \mathbf{z}\in \bfM \, : \  \langle \tilde{\bfx}-\bfx,\mathbf{z}-\bfx\rangle \leq0 \right)\,.
\end{equation}
$P_\bfM$ is nonexpansive and monotone, that is, for all $\bfx,\tilde{\bfx}\in\bfX$,
\begin{equation}\label{Pnonexp}
\|P_\bfM(\bfx)-P_\bfM(\tilde{\bfx})\|\leq \|\bfx-\tilde{\bfx}\|
\end{equation}
and
\begin{equation}\label{Pmon}
\langle P_\bfM(\bfx)-P_\bfM(\tilde{\bfx}),\bfx-\tilde{\bfx}\rangle\geq\|P_\bfM(\bfx)-P_\bfM(\tilde{\bfx})\|^2
\end{equation}
as well as continuous and, in general, nonlinear.

Discretization by restriction to a linear subspace $\bfX_N\subseteq\bfX$ can be easily done by replacing $\bfF:\bfX\to\bfY$ by its restriction \
\begin{align*}
\bfF_N:=\bfF\vert_{\bfX_N}:\bfX_N\to\bfY.
\end{align*}
In our case, $\bfX_N=(Z\times\mathcal{V})^K\times\mathcal{C}_N$, so $\bfX_N$ is a linear space in case of a linear activation function $\sigma$
 albeit not necessarily finite dimensional (for approximation on manifolds, see e.g. \cite{Alberti}). It yields the $k$-th iterate 
\begin{equation}\label{LWprojdiscr}
\bfx_{N,k+1}^\delta=P_\bfM \Bigl(\bfx_{N,k}-\bfF_N'(\bfx_{N,k}^\delta)^\Hadj(\bfF_N(\bfx_{N,k}^\delta)-\bfy^\delta)\Bigr)
\end{equation}
in $\bfX_N$. Doing so, we use the Hilbert space adjoint of $\bfF_N'(\bfx_{N,k}^\delta):\bfX_N\to\bfY$ that is uniquely determined by the identity 
\[
\langle \bfF_N'(\bfx_{N,k}^\delta)^\Hadj \bfy,\bfx_N\rangle = \langle \bfy, \bfF_N'(\bfx_{N,k}^\delta)\bfx_N\rangle \quad
\mbox{ for all }\bfy\in\bfY, \bfx_N\in\bfX_N.
\]
Therefore, the adjoint in the discretized and projected Landweber \eqref{LWprojdiscr} equals to
\begin{align}
\bfF_N'(\bfx_{N,k}^\delta)^\Hadj:\bfY\to\bfX_N \qquad \bfF_N'(\bfx_{N,k}^\delta)^\Hadj=P_{\bfX_N} \bfF'(\bfx_{N,k}^\delta)^\Hadj,
\end{align}
the concatenation of $\bfF'(\bfx_{N,k}^\delta)^\Hadj$ with the orthogonal projection $P_{\bfX_N}:\bfX\to\bfX_N$ onto $\bfX_N$ in the Hilbert space $\bfX$.

\subsection{Convergence} 
Also for the discretized and projected Landweber iteration, we will show that with an appropriate choice of the stopping index $k_*(\delta)$ and discretization parameter $N(\delta)$ it is a regularization method.

We denote by $\bfx^\dagger\in\bfX$ a solution of the inverse problem with exact data, that is, $\bfF(\bfx^\dagger)=\bfy$, by $\bfx_{\infty,k}$ the iterates in $\bfX$ according to \eqref{LWproj}, and make the following assumptions
\begin{assumption}\label{ass3}
\begin{enumerate}[label=(L\arabic*)]
\item\label{LW_a} Approximation by $\bfX_N$: 
There exists a sequence $(\bfx_N^\dagger)_{N\in\mathbb{N}}$, $\bfx_N^\dagger\in \bfX_N\cap \bfM$ such that
for some $\bar{d}$ 
\[
\begin{aligned}
&d_N:=\|\bfx_N^\dagger-\bfx^\dagger\|\leq\bar{d}
\,, \quad
m_N:=\|\bfF_N(\bfx_N^\dagger)-\bfy\|=\|\bfF(\bfx_N^\dagger)-\bfF(\bfx^\dagger)\|\to0\,, \\
&\mbox{ and }s_N:=\sup_{\bfx\in B_R(\bfx^\dagger)}\|(I-P_{\bfX_{N}})\bfF'(\bfx)^\Hadj(\bfF(\bfx)-\bfy)\|\to0 \mbox{ as }
\mbox{ as }N\to\infty;
\end{aligned}
\]
\item\label{LW_b} Convergence and boundedness of the starting values 
\[
\begin{aligned}
&\|\bfx_{N,0}^\delta-\bfx_{\infty,0}\|\leq\rho_N\to0\mbox{ as }N\to\infty\,,\\
&\|\bfx_{N,0}^\delta-\bfx^\dagger\|\leq\rho \,, \quad
\|\bfF(\bfx_{N,0}^\delta) -\bfy\|\leq \tilde{\rho}\,,
\quad \mbox{ for all }N\in\mathbb{N}
\end{aligned}
\]
(e.g., 
by setting $\bfx_{N,0}^\delta:=P_{\bfX_N}\bfx_0$);
\item\label{LW_c} Local boundedness and tangential cone condition on $\bfF$ as well as Lipschitz continuity of $\bfF'$:
There exists $R>\rho+2\bar{d}$, $\mu_R>0$, $M_R>0$, $K_R>0, L_R>0$, such that for all $x\in B_R(\bfx^\dagger)$ and for all $N\in\mathbb{N}$
\begin{equation}\label{MR}
\|\bfF'(\bfx)\|\leq M_R\leq\sqrt{2}
\end{equation}
as well as
\begin{equation}\label{tcc1}
\begin{aligned}
&2\langle \bfF(\bfx)-\bfF(\bfx_N^\dagger),\bfF'(\bfx)(\bfx-\bfx_N^\dagger)\rangle
\geq (M_R^2+\mu_R)\|\bfF(\bfx)-\bfF(\bfx_N^\dagger)\|^2\\
&\|\bfF'(\bfx)(\bfx-\bfx_N^\dagger)\|
\leq K_R\|\bfF(\bfx)-\bfF(\bfx_N^\dagger)\|\,;
\end{aligned}
\end{equation}
\begin{equation}\label{LR}
\|\bfF'(\bfx)-\bfF'(\tilde{\bfx})\|\\
\leq L_R\|\bfx-\tilde{\bfx}\|\,;
\end{equation}
\item\label{LW_d} Asymptotics of the parameters as $\delta\to0$:
\[
\begin{aligned}
&k_*(\delta)\to\infty\,, \quad N(\delta)\to\infty\\
&(\tfrac{4}{\mu_R}K_R+(1+\tfrac{4}{\mu_R}M_R^2)M_R^2)\, k_*(\delta)\,  (m_{N(\delta)}+\delta)^2\leq (R-\bar{d})^2-(\rho+\bar{d})^2; 
\end{aligned}
\]
\item\label{LW_e} The mapping $\bfx\mapsto P_{\bfM}\Bigl(\bfx-\bfF'(\bfx)^\Hadj(\bfF(\bfx)-\bfy)\Bigr)-\bfx$  is weakly sequentially closed;
\item\label{LW_f} For all $y\bf\in B_{\bar{\delta}}(\bfy)$ the mapping $x\mapsto\bfF'(\bfx)^\Hadj(\bfF(\bfx)-\bfz)$ is Lipschitz continuous with constant $L<2$.
\end{enumerate}
\end{assumption}
\begin{remark}[On \ref{LW_a}, \blue{approximation by NNs}]
Note that the first part of \ref{LW_a} (boundedness of $d_N$ and convergence of $m_N$) only requires approximation of the single element $\bfx^\dagger$. By smoothness assumptions on $\bfx^\dagger$, this assumption therefore can be achieved, even with rates for discretization by NNs, \blue{as mentioned in Discussion \ref{appox-rate}}. 

The second part of \ref{LW_a} which is supposed to hold for all $\bfx\in B_R(\bfx^\dagger)$ can be obtained by using the fact that $\bfF'(\bfx)^\Hadj$ is a smoothing operator and therefore even norm convergence $\|(I-P_{\bfX_{N}})\bfF'(\bfx)^\Hadj\|\to0$ follows from error estimates of $I-P_{\bfX_{N}}$ under a priori regularity conditions.
\end{remark}
\begin{remark}[on \ref{LW_c}, tangential cone condition]
A sufficient condition for \eqref{tcc1} is the classical tangential cone condition (cf. \cite{KalNeuSch08})
\begin{equation}\label{tcc}
\begin{aligned}
&\|\bfF(\bfx)-\bfF(\bfx_N^\dagger)-\bfF'(\bfx)(\bfx-\bfx_N^\dagger)\|
\leq c_{tc}\|\bfF(\bfx)-\bfF(\bfx_N^\dagger)\|\,.
\end{aligned}
\end{equation}
for some $c_{tc}<1$ independent of $x\in B_R(\bfx^\dagger)$ with and all $N\in\mathbb{N}$, since by the inverse triangle inequality it is readily checked that we can then set $K_R=1+c_{tc}$ and $M_R+\mu_r=
1-c_{tc}^2+(1-c_{tc})^2$.
\end{remark}

We start with an estimate on the propagated noise and discretization error.
\begin{lemma}\label{lem:noiseprop}
Under conditions \ref{LW_a}, \eqref{MR}, \eqref{LR}, for any $k\in\mathbb{N}$ the estimates
\[
\|\bfx_{N,k}^\delta-\bfx_{\infty,k}\|
\leq (1+\tfrac{5}{2}M_RL_R R)^k\rho_N
+ \frac{2}{5M_RL_R R} (1+\tfrac{5}{2}M_RL_R R)^k (s_N + M_R\delta)
\]
and 
\[
\begin{aligned}
&\|\bfF_N'(\bfx_{N,k}^\delta)^\Hadj(\bfF_N(\bfx_{N,k}^\delta)-\bfy^\delta)
-\bfF'(\bfx_{\infty,k})^\Hadj(\bfF(\bfx_{\infty,k})-\bfy)\|\\
&\leq s_N+M_R(M_R +\tfrac32 L_R R)\|\bfx_{N,k}^\delta-\bfx_{\infty,k}\| + M_R\delta
\end{aligned}
\]
hold, provided that for all $\ell\leq k$, $\bfx_{N,\ell}^\delta$, $\bfx_{\infty,\ell}\in B_R(\bfx^\dagger)$.
\end{lemma}

{\it Proof.}
We make use of the recursions 
\[
\begin{aligned}
\bfx_{N,k+1}^\delta-\bfx^\dagger 
&=\Bigl(I-P_{\bfX_N}{K_{N,k}^\delta}^\Hadj\bar{K}_{N,k}^\delta\Bigr)(\bfx_{N,k}^\delta-\bfx^\dagger)
+P_{\bfX_N}{K_{N,k}^\delta}^\Hadj(\bfy^\delta-\bfy)\\
\bfx_{\infty,k+1}-\bfx^\dagger 
&=\Bigl(I-{K_{\infty,k}}^\Hadj\bar{K}_{\infty,k}\Bigr)(\bfx_{\infty,k}-\bfx^\dagger)
\end{aligned}
\]
where $K_{N,k}^\delta=\bfF'(\bfx_{N,k}^\delta)$, 
$\bar{K}_{N,k}^\delta=\int_0^1\bfF'(\bfx^\dagger+\theta(\bfx_{N,k}^\delta-\bfx^\dagger))\, d\theta$, 
$K_{\infty,k}=\bfF'(\bfx_{\infty,k})$,
$\bar{K}_{\infty,k}=\int_0^1\bfF'(\bfx^\dagger+\theta(\bfx_{\infty,k}-\bfx^\dagger))\, d\theta$. 
This yields
\[
\begin{aligned}
&\bfx_{N,k+1}^\delta-\bfx_{\infty,k+1} \\
&=\Bigl(I-{K_{\infty,k}}^\Hadj K_{\infty,k}\Bigr)(\bfx_{N,k}^\delta-\bfx_{\infty,k})
+{K_{\infty,k}}^\Hadj(K_{\infty,k}-\bar{K}_{\infty,k})(\bfx_{N,k}^\delta-\bfx_{\infty,k})\\
&\quad+\Bigl({K_{\infty,k}}^\Hadj(\bar{K}_{\infty,k}-\bar{K}_{N,k}^\delta)
+({K_{\infty,k}}-{K_{N,k}^\delta})^\Hadj\bar{K}_{N,k}^\delta\Bigr)(\bfx_{N,k}^\delta-\bfx^\dagger)\\
&\quad+(I-P_{\bfX_N}){K_{N,k}^\delta}^\Hadj\bar{K}_{N,k}^\delta(\bfx_{N,k}^\delta-\bfx^\dagger)
+P_{\bfX_N}{K_{N,k}^\delta}^\Hadj(\bfy^\delta-\bfy) 
\end{aligned}
\]
thus by $M_R\leq\sqrt{2}$, which implies $\|I-{K_{\infty,k}}^\Hadj K_{\infty,k}\|\leq1$
\[
\begin{aligned}
\|\bfx_{N,k+1}^\delta-\bfx_{\infty,k+1}\| 
&\leq (1+\tfrac{5}{2}M_RL_R R)\|\bfx_{N,k}^\delta-\bfx_{\infty,k}\|
+ s_N + M_R\delta\\
&\leq (1+\tfrac{5}{2}M_RL_R R)^{k+1}\|\bfx_{N,0}^\delta-\bfx_{\infty,0}\|
+ \sum_{j=0}^k (1+\tfrac{5}{2}M_RL_R R)^j (s_N + M_R\delta)
\end{aligned}
\]
Moreover,
\[
\begin{aligned}
&\bfF_N'(\bfx_{N,k}^\delta)^\Hadj(\bfF_N(\bfx_{N,k}^\delta)-\bfy^\delta)
-\bfF'(\bfx_{\infty,k})^\Hadj(\bfF(\bfx_{\infty,k})-\bfy)\\
&=P_{\bfX_N}{K_{N,k}^\delta}^\Hadj\bar{K}_{N,k}^\delta(\bfx_{N,k}^\delta-\bfx^\dagger) 
-{K_{\infty,k}}^\Hadj\bar{K}_{\infty,k}(\bfx_{\infty,k}-\bfx^\dagger)+{P_{\bfX_N}{K_{N,k}^\delta}^\Hadj(\bfy-\bfy^\delta)}\\
&=-(I-P_{\bfX_N}){K_{N,k}^\delta}^\Hadj\bar{K}_{N,k}^\delta(\bfx_{N,k}^\delta-\bfx^\dagger)
+{K_{N,k}^\delta}^\Hadj\bar{K}_{N,k}^\delta(\bfx_{N,k}^\delta-\bfx_{\infty,k})
+{P_{\bfX_N}}{K_{N,k}^\delta}^\Hadj(\bfy-\bfy^\delta)\\
&\quad
+\Bigl({K_{N,k}^\delta}^\Hadj(\bar{K}_{N,k}^\delta)-\bar{K}_{\infty,k})
+(K_{N,k}^\delta-{K_{\infty,k}})^\Hadj \bar{K}_{\infty,k}\Bigr)
(\bfx_{\infty,k}-\bfx^\dagger)
\end{aligned}
\]
\begin{flushright}
$\diamondsuit$
\end{flushright}

\begin{remark}
In the linear case $\bfF(\bfx)=K\bfx$ with $\|K\|\leq1$ the much better estimates
\[
\|\bfx_{N,k}^\delta-\bfx_{\infty,k}\|
\leq \|\bfx_{N,0}^\delta-\bfx_0\|+k\bigl(\|(I-P_{\bfX_N})K^*K\|R+\delta\bigr)
\]
and 
\[
\begin{aligned}
&\|\bfF_N'(\bfx_{N,k}^\delta)^\Hadj(\bfF_N(\bfx_{N,k}^\delta)-\bfy^\delta)
-\bfF'(\bfx_{\infty,k})^\Hadj(\bfF(\bfx_{\infty,k})-\bfy)\|\\
&\leq \frac{1}{k+1}\|\bfx_{N,0}^\delta-\bfx_0\|+\Bigl(1+\sum_{j=0}^{k-1}\frac{1}{j+1}\Bigr)\bigl(\|(I-P_{\bfX_N})K^*K\|R+\delta\bigr)
\end{aligned}
\]
can be easily verified by means of spectral theoretic methods. 
More precisely, we use the fact that $\|K^*K(I-K^*K)^j\|\leq\frac{1}{j+1}$, and the identities 
\[\begin{aligned}
&\bfx_{N,k}^\delta-\bfx_{\infty,k}
=(\bfx_{N,k}^\delta-\bfx^\dagger)-(\bfx_{\infty,k}-x^\dagger)\\
&=(I-K^*K)^k(\bfx_{N,0}^\delta-x^\dagger)- (I-K^*K)^k(\bfx_{\infty,0}-\bfx^\dagger)\\
&\quad+\sum_{j=0}^{k-1}(I-K^*K)^{j}
\Bigl((I-P_{\bfX_N})K^*K(\bfx_{N,k-j-1}^\delta-x^\dagger)+P_{\bfX_N}K^*(\bfy^\delta-\bfy)\Bigr)
\end{aligned}\]
\[\begin{aligned}
&P_{\bfX_N}K^*(K\bfx_{N,k}^\delta-\bfy^\delta)-K^*(\bfx_{\infty,k}-\bfy)\\
&=K^*K(\bfx_{N,k}^\delta-\bfx_{\infty,k})-(I-P_{\bfX_N})K^*K(\bfx_{N,k}^\delta-x^\dagger)-P_{\bfX_N}K^*(\bfy^\delta-\bfy)\\
&=K^*K(I-K^*K)^k(\bfx_{N,0}^\delta-\bfx_{\infty,0})\\
&\quad
+\sum_{j=0}^{k-1}K^*K(I-K^*K)^{j}
\Bigl((I-P_{\bfX_N})K^*K(\bfx_{N,k-j-1}^\delta-x^\dagger)+P_{\bfX_N}K^*(\bfy^\delta-\bfy)\Bigr)\\
&\quad -(I-P_{\bfX_N})K^*K(\bfx_{N,k}^\delta-x^\dagger)-P_{\bfX_N}K^*(\bfy^\delta-\bfy)
\end{aligned}\]

They can actually be transferred to the nonlinear setting under an adjoint range invariance condition on $\bfF$, which is a stronger assumption than the tangential cone condition, similarly to the convergence rates estimates in \cite{HNS95}.
However, in our example, this \blue{assumption} does not seem to be verifiable, whereas the tangential cone condition can be established, see below.
\end{remark}

While uniform boundedness of the iterates can be shown under the assumptions \ref{LW_a}-\ref{LW_d}, in order to control the propagated noise in the iterates, we will therefore have to additionally impose
\begin{assumption}
\begin{equation}\label{klogdelta_delta0}
(1+\tfrac{5}{2}M_RL_R R)^k\rho_{N(k)}\to0\,, \quad
(1+\tfrac{5}{2}M_RL_R R)^{k} s_{N(k)}  \to0\mbox{ as }k\to\infty
\end{equation}
in case of exact data $\delta=0$ and 
\begin{equation}\label{klogdelta}
(1+\tfrac{5}{2}M_RL_R R)^{k_*(\delta)}\rho_{N(\delta)}\to0\,, \quad
(1+\tfrac{5}{2}M_RL_R R)^{k_*(\delta)} (s_{N(\delta)} + M_R\delta) \to0\mbox{ as }\delta\to0\,.
\end{equation}
\end{assumption}

\begin{proposition}\label{prop:LW}
Under the above assumptions \ref{LW_a}-\ref{LW_f} with $\bfM$ closed and convex, the iterates are well-defined by \eqref{LWprojdiscr} and remain in $B_R(\bfx^\dagger)$. 

Under the additional condition \eqref{klogdelta}, we also have weak subsequential convergence of
$\bfx_{N(\delta),k_*(\delta)}^\delta$ to a solution of \eqref{Fxy} as $\delta\to0$.

With exact data $\delta=0$ and $N=N(k)$ chosen according to \eqref{klogdelta_delta0}, we have weak subsequential convergence of $\bfx_{N(k),k}$ to a solution of \eqref{Fxy} as $k\to\infty$.
\end{proposition}
{\it Proof.}
We follow the classical monotonicity proof from \cite{HNS95}, see also \cite{Kindermann2017}, but do so with $F_N$ instead of $F$ so that we can exploit the identity $\langle \bfF_N'(\bfx_{N,k}^\delta)^\Hadj(\bfF_N(\bfx_{N,k}^\delta)-\bfy^\delta),(\bfx_{N,k}^\delta-\bfx_N^\dagger)\rangle
=\langle \bfF_N(\bfx_{N,k}^\delta)-\bfy^\delta,\bfF_N'(\bfx_{N,k}^\delta)(\bfx_{N,k}^\delta-\bfx_N^\dagger)\rangle$ in the first equality below. It is also for this reason that we had to introduce the auxiliary variable $x_N^\dagger$ as a substitute for $x^\dagger$ in $\bfX_N$.
Therewith we obtain, for arbitrary $N\in\mathbb{N}$, using the fact that we can skip the subscript $N$ when applying $\bfF_N$ to an element of $\bfX_N$ and nonexpansivity \eqref{Pnonexp} together with the fact that $P_\bfM(\bfx_N^\dagger)=\bfx_N^\dagger$
\[
\begin{aligned}
&\|\bfx_{N,k+1}^\delta-\bfx_N^\dagger\|^2-\|\bfx_{N,k}^\delta-\bfx_N^\dagger\|^2\\
&=\|\bfF_N'(\bfx_{N,k}^\delta)^\Hadj(\bfF_N(\bfx_{N,k}^\delta)-\bfy^\delta)\|^2
-2\langle \bfF_N(\bfx_{N,k}^\delta)-\bfy^\delta,\bfF_N'(\bfx_{N,k}^\delta)(\bfx_{N,k}^\delta-\bfx_N^\dagger)\rangle\\
&=\|\bfF_N'(\bfx_{N,k}^\delta)^\Hadj(\bfF_N(\bfx_{N,k}^\delta)-\bfF_N(\bfx_N^\dagger))\|^2
-2\langle \bfF_N(\bfx_{N,k}^\delta)-\bfF_N(\bfx_N^\dagger),\bfF_N'(\bfx_{N,k}^\delta)(\bfx_{N,k}^\delta-\bfx_N^\dagger)\rangle\\
&\qquad
+2\langle \bfy^\delta-\bfF_N(\bfx_N^\dagger),\bfF_N'(\bfx_{N,k}^\delta)\Bigl((\bfx_{N,k}^\delta-\bfx_N^\dagger)-\bfF_N'(\bfx_{N,k}^\delta)^\Hadj(\bfF_N(\bfx_{N,k}^\delta)-\bfF_N(\bfx_N^\dagger))\Bigr)\rangle\\
&\qquad
+\|\bfF_N'(\bfx_{N,k}^\delta)^\Hadj(\bfy^\delta-\bfF_N(\bfx_N^\dagger))\|^2\\
&\leq -\mu_R\|\bfF(\bfx_{N,k}^\delta)-\bfF(\bfx_N^\dagger)\|^2
+(\tfrac{2}{\epsilon}+M_R^2)\|\bfy^\delta-\bfF(\bfx_N^\dagger)\|^2\\
&\qquad+\epsilon \|\bfF'(\bfx_{N,k}^\delta)(\bfx_{N,k}^\delta-\bfx_N^\dagger)\|^2
+\epsilon M_R^4 \|\bfF(\bfx_{N,k}^\delta)-\bfF(\bfx_N^\dagger)\|^2\\
&\leq -(\mu_R-\epsilon K_R-\epsilon M_R^4)\|\bfF(\bfx_{N,k}^\delta)-\bfF(\bfx_N^\dagger)\|^2
+(\tfrac{2}{\epsilon}+M_R^2)(\delta+m_N)^2\\
&\leq-\tfrac{\mu_R}{2} \|\bfF(\bfx_{N,k}^\delta)-\bfF(\bfx_N^\dagger)\|^2
+(\tfrac{4}{\mu_R}K_R+(1+\tfrac{4}{\mu_R}M_R^2)M_R^2)(\delta+m_N)^2
\end{aligned}
\]
provided $\bfx_{N,k}^\delta\in B_R(\bfx^\dagger)$.
Here we have employed Young's inequality in the form 
$2a(b+c)\leq\frac{2}{\epsilon} a^2 + \frac{\epsilon}{2}(b+c)^2 \leq \frac{2}{\epsilon} a^2 + \epsilon b^2 + \epsilon c^2$ with $\epsilon=\frac{\mu_R}{2(K_R+M_R^4)}$.

Summing up for $k$ from zero to $\tilde{k}-1$ we obtain
\begin{equation}\label{sumLW}
\begin{aligned}
&\tfrac{\mu_R}{2} \sum_{k=0}^{\tilde{k}-1} \|\bfF(\bfx_{N,k}^\delta)-\bfF(\bfx_N^\dagger)\|^2
+\|\bfx_{N,\tilde{k}}^\delta-\bfx^\dagger_N\|^2\\
&\leq 
\|\bfx_{N,0}^\delta-\bfx^\dagger_N\|^2 +\tilde{k}(\tfrac{4}{\mu_R}K_R+(1+\tfrac{4}{\mu_R}M_R^2)M_R^2)(\delta+m_N)^2\,,
\end{aligned}
\end{equation}
which by \ref{LW_a}, \ref{LW_b} and \ref{LW_d} inductively implies that the iterates $\bfx_{N(\delta),\tilde{k}}^\delta$ remain in $B_R(\bfx^\dagger)$ for all $\tilde{k}\leq k_*(\delta)$.
Thus $(\bfx_{N(\delta),k_*(\delta)}^\delta)_{\delta>0}$ has a weakly convergent subsequence
\begin{equation}\label{weakconvx}
\bfx_{N(\delta^j),k_*(\delta^j)}^{\delta^j}\rightharpoonup \bar{\bfx}\,.
\end{equation}
and in case $\delta=0$, 
with $N=N(k)$ in place of $N=N(\delta)$
\begin{equation}\label{weakconvxdelta0}
\bfx_{N(k_j),k_j}\rightharpoonup \bar{\bfx}\,.
\end{equation}
Since $\bfM$ is closed and convex, hence weakly closed, $\bar{\bfx}$ is contained in $\bfM$.

\medskip

To prove that this limit solves the inverse problem, like in \cite[Lemma 3.1]{Kindermann2017} with 
\[
J_N^\delta(\bfx)=\frac12\|\bfF_N(\bfx)-\bfy^\delta\|^2\,, \qquad 
\Delta_{N,k}^\delta:=\bfx_{N,k+1}^\delta-\bfx_{N,k}^\delta=
P_\bfM\Bigl(\bfx_{N,k}^\delta-{J_N^\delta}'(\bfx_{N,k})\Bigr)-\bfx_{N,k}^\delta
\] 
and \ref{LW_f}, which implies that for all $\bfx_N, \tilde{\bfx}_N\in B_R(\bfx^\dagger)$
\[
\begin{aligned}
&\|{J_N^\delta}'(\bfx_N)-{J_N^\delta}'(\tilde{\bfx}_N)\|
= \|\bfF_N'(\bfx_N)^\Hadj(\bfF_N(\bfx_N)-\bfy^\delta)-\bfF_N'(\tilde{\bfx}_N)^\Hadj(\bfF_N(\tilde{\bfx}_N)-\bfy^\delta)\|\\
&= \|P_{\bfX_N}\Bigl(\bfF'(\bfx_N)^\Hadj(\bfF(\bfx_N)-\bfy^\delta)-\bfF'(\tilde{\bfx}_N)^\Hadj(\bfF(\tilde{\bfx}_N)-\bfy^\delta)\Bigr)\|
\leq L\|\bfx_N-\tilde{\bfx}_N\|
\end{aligned}.
\]
We then obtain, for any $k$ \blue{and for both the discretized problem in $\bfX_N$ as well as non-discretized problem in $\bfX$ (i.e. $N=\infty$)},
\[
\begin{aligned}
&{J_N^\delta}(\bfx_{N,k+1})-{J_N^\delta}(\bfx_{N,k})
=\int_0^1 {J_N^\delta}'(\bfx_{N,k}^\delta+\theta\Delta_{N,k}^\delta)\Delta_{N,k}^\delta\, d\theta\\
&\leq-\|\Delta_{N,k}^\delta\|^2 + \int_0^1 ({J_N^\delta}'(\bfx_{N,k}^\delta+\theta\Delta_{N,k}^\delta)-{J_N^\delta}'(\bfx_{N,k}^\delta))\Delta_{N,k}^\delta\, d\theta
\leq -(1-\tfrac{L}{2}) \|\Delta_{N,k}^\delta\|^2\,,
\end{aligned}
\]
where we have used the fact that monotonicity \eqref{Pmon} 
with $\bfx=\bfx_{N,k}^\delta-{J_N^\delta}'(\bfx_{N,k}^\delta)$, $\tilde{\bfx}=\bfx_{N,k}^\delta=P_\bfM(\bfx_{N,k}^\delta)$ since  $\bfx_{N,k}^\delta\in\bfM$
implies 
\[
-{J_N^\delta}'(\bfx_{N,k}^\delta)\Delta_{N,k}^\delta \geq \|\Delta_{N,k}^\delta\|^2\,.
\]
After summation and by \ref{LW_b} and $J_N^\delta\geq0$ this implies that  
\begin{equation}\label{sumNinfty}
\begin{aligned}
\sup_{\delta\in(0,\bar{\delta}], \, N\in\mathbb{N}}
\sum_{k=0}^{\infty} \|\Delta_{N,k}^\delta\|^2 
\leq\frac{1}{2-L} \sup_{N\in\mathbb{N}}\|\bfF(\bfx_{N,0}^\delta)-\bfy^\delta\|^2
\leq \frac{1}{2-L} (\tilde{\rho}+\bar{\delta})^2\,,
\end{aligned}
\end{equation}
where 
\[
\begin{aligned}
\Delta_{N,k}^\delta
&=P_{\bfM}\Bigl(\bfx_{N,k}-\bfF_{N}'(\bfx_{N,k}^\delta)^\Hadj(\bfF_{N}(\bfx_{N,k}^\delta)-\bfy^\delta)\Bigr) - \bfx_{N,k}\\
&=P_{\bfM}\Bigl(\bfx_{N,k}-P_{\bfX_{N}}\bfF'(\bfx_{N,k})^\Hadj(\bfF(\bfx_{N,k})-\bfy^\delta)\Bigr)
- \bfx_{N,k}
\end{aligned}
\]

In particular, in case $\delta=0$ (thus skipping the superscript $\delta$ 
and setting $N=N(k)$) 
with nonexpansivity of $P_{\bfM}$ 
\[
\begin{aligned}
\|\Delta_{N,k}^\delta\|&=\|P_{\bfM}\Bigl(\bfx_{N(k),k}-\bfF'(\bfx_{N(k),k})^\Hadj(\bfF(\bfx_{N(k),k})-\bfy)\Bigr) - \bfx_{N(k),k}\|\\
&=\|P_{\bfM}\Bigl(\bfx_{N(k),k}-P_{\bfX_N}\bfF'(\bfx_{N(k),k})^\Hadj(\bfF(\bfx_{N(k),k})-\bfy)\Bigr) - \bfx_{N(k),k} \\
&\qquad+
P_{\bfM}\Bigl(\bfx_{N(k),k}-\bfF'(\bfx_{N(k),k})^\Hadj(\bfF(\bfx_{N(k),k})-\bfy)\Bigr) \\
&\qquad-P_{\bfM}\Bigl(\bfx_{N(k),k}-P_{\bfX_N}\bfF'(\bfx_{N(k),k})^\Hadj(\bfF(\bfx_{N(k),k})-\bfy)\Bigr)
\|\\
&\leq \|\Delta_{N(k),k}\|+ \|(I-P_{\bfX_{N(k)}})\bfF'(\bfx_{N(k),k})^\Hadj(\bfF(\bfx_{N(k),k})-\bfy)\|\\
&\leq \|\Delta_{\infty,k}\|+ 2\|(I-P_{\bfX_{N(k)}})\bfF'(\bfx_{N(k),k})^\Hadj(\bfF(\bfx_{N(k),k})-\bfy)\| +2\|\bfx_{N(k),k}-\bfx_{\infty,k}\|\\
&\qquad +
\|\bfF'(\bfx_{N(k),k})^\Hadj(\bfF(\bfx_{N(k),k})-\bfy)-\bfF'(\bfx_{\infty,k})^\Hadj(\bfF(\bfx_{\infty,k})-\bfy)\|\\
&\to 0\mbox{ as }k\to\infty\,,
\end{aligned}
\]
due to \eqref{sumNinfty}, \ref{LW_a}
and \eqref{klogdelta_delta0}, according to Lemma~\ref{lem:noiseprop}.

Thus from \eqref{weakconvxdelta0}, $\bar{x}\in\bfM$ and \ref{LW_e} we get 
$P_{\bfM}\Bigl(\bar{\bfx}-\bfF'(\bar{\bfx})^\Hadj(\bfF(\bar{\bfx})-\bfy)\Bigr) - \bar{\bfx}=0$, 
hence due to \eqref{proj} with 
$\tilde{\bfx}=\bar{\bfx}-\bfF'(\bar{\bfx})^\Hadj(\bfF(\bar{\bfx})-\bfy)$,
$\bfx=\bar{\bfx}$, $\mathbf{z}=\bfx^\dagger$
and \ref{LW_c}
\[
0\geq\langle \bfF'(\bar{\bfx})(\bar{\bfx}-\bfx^\dagger),\bfF(\bar{\bfx})-\bfF(\bfx^\dagger)\rangle
\geq\frac{M_R^2+\mu_R}{2}\|\bfF(\bar{\bfx})-\bfF(\bfx^\dagger)\|^2\,.
\]
This gives subsequential convergence of $\bfx_{N(\delta),k(\delta)}$ to a solution $x^\dagger$ of \eqref{Fxy} as $k\to\infty$ with exact data, \blue{for both the discretized and the nondiscretized problem.}

Convergence with noisy data can be concluded from Lemma~\ref{lem:noiseprop} under the more restrictive assumption \eqref{klogdelta}. \blue{Indeed, in the decomposition
\begin{align*}
\bfx_{N(\delta),k(\delta)}^\delta-\bfx^\dagger=(\bfx_{N(\delta),k(\delta)}^\delta-\bfx_{\infty,k(\delta)})+(\bfx_{\infty,k(\delta)}-\bfx^\dagger),
\end{align*}
convergence of the first term follows from Lemma~\ref{lem:noiseprop} and the rule \eqref{klogdelta}, while weak convergence of the second term is a consequence of the result with exact data in the nondiscretized setting $N=\infty$, that we have just proven above.}
\begin{flushright}
$\diamondsuit$
\end{flushright}

\subsection{Adjoint}\label{sec:Land-adjoint} We will now write out $\bfF'(\bfx)^\Hadj$ and the expression $\bfF'(\bfx)^\Hadj(\bfF(\bfx)-\bfy)$ that plays a role both in the definition of the Landweber iteration and in the verification of the conditions \ref{LW_e}, \ref{LW_f}.

To do so, we recall the setting
{\allowdisplaybreaks
\begin{align}\label{VWY_LW}
&\bfF(\bfx)=\left(\begin{array}{c}\dot{u}-F(\lambda,u)-f(\alpha,u)\\u(0)-u_0\\Mu\end{array}\right)\in \Vc\times H\times\Yc \,, \quad 
\bfy =\left(\begin{array}{c}0\\ 0\\ y\end{array}\right)\,, \quad 
\tilde{\bfy} =\left(\begin{array}{c}\wtil\\ \htil\\ \ytil\end{array}\right)\,,\nonumber\\[1.5ex]
&\bfF'(\bfx)\tilde{\bfx}=\left(\begin{array}{l}-F_\lambda(\lambda,u)\latil-f_\alpha(\alpha,u)\altil+\dot{\util}-F_u(\lambda,u)\util-f_u(\alpha,u)\util-\ftil(\alpha,u)\\\util(0)-\util_0\\M\util\end{array}\right) \,, \quad 
\nonumber\\
&\bfx =\left(\begin{array}{c}\lambda\\u_0\\ \alpha\\ u\\ f\end{array}\right)\,, \quad 
\tilde{\bfx} =\left(\begin{array}{c}\latil\\ \util_0\\ \altil\\ \util\\ \ftil\end{array}\right) \,, \quad 
\bfF'(\bfx)^\Hadj\tilde{\bfy} = \left(\begin{array}{c}\mutil\\ \vtil_0\\ \betil\\ \vtil\\ \gtil\end{array}\right) \in X\times H\times \R^n\times\Vc\times\Cc, \nonumber\\
&\mathcal{V}=H^1(0,T;V^*)\cap L^2(0,T;V)\,, \quad
\mathcal{W}=L^2(0,T;V^*)\,, \quad
\mathcal{Y}=L^2(0,T;Y)\,, \\
&\langle \util,\vtil \rangle_\mathcal{V} 
=\int_0^T \Bigl(\langle \dot{\util}(t),\dot{\vtil}(t)\rangle_{V^*}+\langle\util(t),\vtil(t)\rangle_V\Bigr)\, dt \nonumber\\ 
&\qquad\quad=\int_0^T \Bigl(\langle \dot{\util}(t),I_V\dot{\vtil}(t)\rangle_{V^*,V} +\langle D_V\vtil(t),\util(t)\rangle_{V^*,V}\Bigr)\, dt, \nonumber \\
&\langle w,\wtil \rangle_\mathcal{W} 
=\int_0^T \langle w(t),\wtil(t)\rangle_{V^*}\, dt 
=\int_0^T \langle w(t),I_V\wtil(t)\rangle_{V^*,V}\, dt  \nonumber
\end{align}}
with the Riesz isomorphisms $I_V:V^*\to V$, $D_V:V\to V^*$ and $V\embed H\embed V^*$ forming a Gelfand triple, and  a Hilbert parameter space $X$. We use the integration by parts identity
\[
\int_0^T \Bigl(\langle \dot{\util}(t),z(t)\rangle_{V^*,V}+\langle \dot{z}(t),\util(t)\rangle_{V^*,V}
=\langle \util(T),z(T)\rangle_H -\langle \util(0),z(0)\rangle_H\,.
\]
Moreover, in order to work in a Hilbert space setting, we will use the Bochner Sobolev space
\begin{equation}\label{C_LW}
\mathcal{C}=H^\ell(\mathbb{R}^n,H^r(\mathbb{R}))\,, \quad
\langle\ftil,\gtil\rangle_{\mathcal{C}}:=\int_{\mathbb{R}^n} (1+|\kappa|^2)^\ell\int_{\mathbb{R}}
(1+|\omega|^2)^r(\mathcal{F}\ftil)(\kappa,\omega)\overline{\mathcal{F}\gtil)(\kappa,\omega)}\, d\omega\, d\kappa
\end{equation}
with $\ell$, $r$ large enough to allow for $\mathcal{C}\subseteq C(\mathbb{R}^n,\mathbb{R})\cap C(\mathbb{R}^n;W^{1,\infty}(\mathbb{R}))$, see \eqref{tangconef} below, where $\mathcal{F}$ denotes the Fourier transform.
Therewith, the defining identity for the Hilbert space adjoint $\bfF'(\bfx)^\Hadj\tilde{\bfy}$, that is,
\[
0=\langle \bfF'(\bfx)^\Hadj \tilde{\bfy},\tilde{\bfx}\rangle - \langle \tilde{\bfy}, \bfF'(\bfx)\tilde{\bfx}\rangle \quad
\mbox{ for all }\bfx\in\bfX
\]
reads as follows:
{\allowdisplaybreaks
\begin{align*}
0=&\langle\latil,\mutil\rangle_X 
+ \langle\util_0,\vtil_0\rangle_H 
+ \langle\altil,\betil\rangle_{\mathbb{R}^n}
+ \int_0^T \Bigl(\langle \dot{\util}(t),I_V\dot{\vtil}(t)\rangle_{V^*,V} +\langle D_V\vtil(t),\util(t)\rangle_{V^*,V}
\Bigr)\, dt
+ \langle\ftil,\gtil\rangle_{\mathcal{C}}
\\
&+\int_0^T \langle F_\lambda(\lambda,u)\latil+f_\alpha(\alpha,u)\altil-\dot{\util}+F_u(\lambda,u)\util+f_u(\alpha,u)\util+\ftil(\alpha,u),I_V\wtil\rangle_{V^*,V}\,dt\\
&-\langle\util(0)-\util_0,\htil\rangle_H
-\int_0^T\langle M\util,\ytil\rangle_Y\, dt\\
=&\int_0^T \Bigl(\langle -I_V\ddot{\vtil}+I_V\dot{\wtil} +F_u(\lambda,u)^*I_V\wtil +f_u(\alpha,u)^*I_V\wtil -M^\Hadj\ytil +D_V\vtil,\util\rangle_{V^*,V}\,dt\\
&+\langle I_V\dot{\vtil}(T)-I_V\wtil(T),\util(T)\rangle_H
-\langle I_V\dot{\vtil}(0)-I_V\wtil(0)-\htil,\util(0)\rangle_H
+\langle \htil+\vtil_0,\util_0\rangle_H\\
&+\langle \int_0^T F_\lambda(\lambda,u)^*I_V\wtil\, dt+\mutil,\latil\rangle_X
+\langle \int_0^T f_\alpha(\alpha,u)^*I_V\wtil\, dt+\betil,\altil\rangle_{\mathbb{R}^n}\\
&+\int_{\mathbb{R}^n}\int_{\mathbb{R}} 
\Bigl((1+|\kappa|^2)^\ell(1+|\omega|^2)^r \overline{\mathcal{F}\gtil(\kappa,\omega)}\\
&\qquad\qquad+\frac{1}{2\pi^{(n+1)/2}}
\Bigl[\int_0^T\int_\Omega e^{i\kappa\cdot\alpha} e^{i\omega u(x,t)} (I_V\wtil)(x,t)\, dx\,dt
\Bigr]
\Bigr)\,\mathcal{F}\ftil(\beta,\omega)\, d\omega\, d\kappa
\end{align*}}
where we have rewritten 
\[
\ftil(\alpha,u(x,t))=\frac{1}{2\pi^{(n+1)/2}}\int_{\mathbb{R}^n}\int_{\mathbb{R}} 
e^{i\kappa\cdot\alpha} e^{i\omega u(x,t)} \mathcal{F}\ftil(\kappa,\omega) \, d\omega\, d\kappa
\]
by the definition of the Fourier transform.
This leads us to defining 
{\allowdisplaybreaks
\begin{align}\label{Land-adjoint}
&\bfF'(\bfx)^\Hadj\tilde{\bfy} = 
\left(\mutil\,\ \vtil_0\,\ \betil\,\ \vtil\,\ \gtil\right)^T \in X\times H\times \R^n\times\Vc\times\Cc, \nonumber\\
&\mutil=-\int_0^T F_\lambda(\lambda,u)^*I_V\wtil\, dt, \nonumber \\
&\vtil_0= -\htil \nonumber\\
&\betil=-\int_0^T f_\alpha(\alpha,u)^*I_V\wtil\, dt,\\
&\vtil= I_V^{-1}\ztil, \nonumber\\
&\gtil=-\frac{1}{2\pi^{(n+1)/2}}\mathcal{F}^{-1}\Bigl[(1+|\kappa|^2)^{-\ell}(1+|\omega|^2)^{-r} 
\Bigl(\int_0^T\int_\Omega e^{-i\kappa\cdot\alpha} e^{-i\omega u(x,t)} (I_V\wtil)(x,t)\, dx\,dt,
\Bigr)\Bigr], \nonumber
\end{align}}
where $\ztil$ solves the two point boundary value problem
\begin{equation}\label{Land-adjoint-wave}
\begin{split}
&\ddot{\ztil}-D_VI_V^{-1}\ztil=I_V\dot{\wtil} +F_u(\lambda,u)^*I_V\wtil +f_u(\alpha,u)^*I_V\wtil -M^\Hadj\ytil\\
&\dot{\ztil}(0)=I_V\wtil(0)+\htil\,, \quad \dot{\ztil}(T)=I_V\wtil(T)
\end{split}
\end{equation}
and $F_u(\lambda,u)^*, f_u(\alpha,u)^*:V\to V^*, F_\lambda(\lambda,u)^*:V\to X, f_\alpha^*(\alpha,u):V\to\R^n$ and $M^*:Y\to V^*$ are Banach space adjoints.
(Note that in case $V=H_0^1(\Omega)$, we have $D_V=I_V^{-1}=-\Delta$ and so the above is a wave equation with the bi-Laplace operator.)

\def\w{\mathbf{w}}
\def\ctil{\tilde{c}}
\def\y{\mathbf{y}}
\def\x{\mathbf{x}}
\def\Fb{\mathbb{F}}
\def\Fc{\mathcal{F}}

\subsection{Discussion of the Assumptions for Application \ref{IP_c-problem}}\label{Land-dis}

We focus on the special case  from the Application \eqref{IP_c-problem}, \eqref{Muy}, that is,
\begin{align}\label{Land-dis-eq}
\bfF(\x)=\bfF(c,\source,u_0,u,f) = \left(\begin{array}{c}\dot{u} -\Delta u + c u +h(u)- f(u)-\source\\u(0)-u_0\\ Mu\end{array}\right)\,,
\end{align}
with 
\begin{equation}\label{spaces_cprob_LW}
H=L^2(\Omega), \quad V=H_0^1(\Omega), \quad X=X_c\times X_\source, \quad 
X_c=L^2(\Omega),\quad X_\source=V^*,
\end{equation} 
(cf. \eqref{VWY_LW}, \eqref{C_LW} for the resulting spaces $\Vc$, $\Wc$, $\Yc$, $\Cc$) and the known nonlinearity $h\in W^{2,\infty}(B)$.

At the end of this section, we will conclude convergence of Landweber iteration and also of Tikhonov regularization for this application from the analysis of the requirements in the following Sections~\ref{sec:tangcone}--\ref{sec:Lipschitz}.

\subsubsection{Tangential cone condition}\label{sec:tangcone}
\[
\begin{aligned}
&\|\bfF(u,f)-\bfF(\tilde{u},\tilde{f})-\bfF'(u,f)(u-\tilde{u},f-\tilde{f})\|_{\mathcal{W}\times H\times\mathcal{Y}} \\
&\leq \|(c-\tilde{c})(u-\tilde{u})\|_{\mathcal{W}}
+\|h(u)-h(\tilde{u})-h'(u)(u-\tilde{u})\|_{\mathcal{W}}\\
&\qquad+\|f(u)-\tilde{f}(\tilde{u})-f'(u)(u-\tilde{u}) - (f-\tilde{f})(u)\|_{\mathcal{W}}\\
&=I+II+III,
\end{aligned}
\]
where 
\[
\begin{aligned}
III&= \|\tilde{f}(u)-\tilde{f}(\tilde{u})-f'(u)(u-\tilde{u})\|_{\mathcal{W}}\\
&= \|\int_0^1 \Bigl({\tilde{f}}'(u+\theta(\tilde{u}-u))-f'(u)\Bigr)\, d\theta(u-\tilde{u})\|_{\mathcal{W}}\\
&= \|\int_0^1 \Bigl({\tilde{f}}'(u+\theta(\tilde{u}-u))-\tilde{f}'(u)\Bigr)\, d\theta\,(u-\tilde{u}) + (f-\tilde{f})'(u)(u-\tilde{u})\|_{\mathcal{W}}\\
&= \|\int_0^1 \int_0^1{\tilde{f}}''(u+s\theta(\tilde{u}-u))\, ds\,\theta d\theta\,(u-\tilde{u})^2 + (f-\tilde{f})'(u)(u-\tilde{u})\|_{\mathcal{W}}\,.
\end{aligned}
\]
So with full observations $Mu=u$ and a choice of spaces 
\begin{equation}\label{spaces_cprob_LW_1}
\mathcal{Y}=L^2(0,T;L^2(\Omega)),\quad
L^p(0,T;L^p(\Omega))\subseteq \mathcal{W} \mbox{ for }p\in\{1,2\},\quad  
\mathcal{C}\subseteq W^{1,\infty}(B)
\end{equation}
with the embedding constant $C_p$, $p\in\{1,2\}$, and  $\mbox{supp}(u)\cup\mbox{supp}(\tilde{u})\subset B$, we obtain
\begin{equation}\label{tangconef}
\begin{aligned}
&\|\bfF(u,f)-\bfF(\tilde{u},\tilde{f})-\bfF'(u,f)(u-\tilde{u},f-\tilde{f})\|_{\mathcal{W}\times H\times\mathcal{Y}}\\
&\leq C_1 \|c-\tilde{c}\|_{L^2} \|u-\tilde{u}\|_{L^1(L^2)}
+\tfrac12 C_1\|h''\|_{L^\infty(B)} \|u-\tilde{u}\|_{L^2(L^2)}^2\\
&\qquad +\tfrac12 C_1\|\tilde{f}''\|_{L^\infty(B)} \|u-\tilde{u}\|_{L^2(L^2)}^2
+ C_2\|(f-\tilde{f})'\|_{L^\infty(B)} \|u-\tilde{u}\|_{L^2(L^2)}\\
&\leq \Bigl(
C_1 \sqrt{T}\|c-\tilde{c}\|_{L^2} 
+\tfrac12 C_1(\|h''\|_{L^\infty(B)}+\|\tilde{f}''\|_{L^\infty(B)}) 
C_{\mathcal{V}\to\mathcal{Y}}\|u-\tilde{u}\|_{\mathcal{V}}\\
&\hspace*{8cm}+ C_2\|f-\tilde{f}\|_{\mathcal{C}}\Bigr) \|Mu-M\tilde{u}\|_{\mathcal{Y}}\\
&\leq c_{tc} \|\bfF(u,f)-\bfF(\tilde{u},\tilde{f})\|_{\mathcal{W}\times H\times\mathcal{Y}}
\end{aligned}
\end{equation}
for all $(\tilde{u},\tilde{f})=(u_N^\dagger,f_N^\dagger)$, $N\in\mathbb{N}$, $(u,f)\in B_R^{\mathcal{V}\times\mathcal{C}}(u^\dagger,f^\dagger)$, provided $R$, and hence $\rho$, are small enough so that  
\[
\sup_{N\in\mathbb{N}} \Bigl(C_1 \sqrt{T}+\tfrac12 C_1((\|h''\|_{L^\infty(B)}+\|{f_N^\dagger}''\|_{L^\infty(B)}) C_{\mathcal{V}\to\mathcal{Y}} 
+ C_2 \Bigr)\, (R+\bar{d})\leq c_{tc} \,.
\]

\subsubsection{Weak sequential closedness of $\x\mapsto-\Fb'(\x)^*(\Fb(\x)-\y)$}
In this section, we study the more general case, namely the application \eqref{Land-dis-eq} with $f=f(\alpha,u)$. We will derive weak closedness via weak continuity.

In the following, we frequently employ the embeddings \cite[Theorems 1.20, 1.21, Lemmas 7.3, 7.7]{Roubicek}, \cite[Chapter 4]{Adams}, \cite[Chapter 11]{Leoni:2009} cf. \eqref{VWY_LW}, \eqref{spaces_cprob_LW}
\begin{align*}
&\Vc\embed C(0 ,T;L^2(\Omega)),\quad \Vc\compt L^2(0,T;L^{6-\epsilon'}(\Omega)), 0<\epsilon'\leq 5, \quad \Wc^*\embed L^2(0,T;L^{6}(\Omega)),\\
&\Cc=H^s(\R^{n+1})\embed C_b(\R^{n+1}), \,\,\qquad s>(n+1)/2\\
&\Cc=H^s(\R^{n+1})\embed W^{2,\infty}(\R^{n+1}),\quad s>(n+1)/2+2
\end{align*}
as well as the H\"older inequalities 
\[
\int_\Omega abc\,dx\leq \|a\|_{L^{3/2}}\|b\|_{L^6}\|c\|_{L^6}, \quad \int_\Omega abc\,dx\leq \|a\|_{L^2}\|b\|_{L^3}\|c\|_{L^6}.
\]

Let $u_n\overset{\Vc}{\rightharpoonup}u, f_n\overset{\Cc}{\rightharpoonup}f, (c,\source)_n\overset{X}{\rightharpoonup}(c,\source), (u_0)_n\overset{H}{\rightharpoonup}u_0, \alpha_n\overset{\R^n}{\to}\alpha$. We first show weak continuity of the model operator $\x\mapsto(\Fb(\x)-y)$. 
\begin{proposition}\label{prop:weakcont}
The operator $\Fb$ defined by \eqref{Land-dis-eq} is weakly continuous on the spaces \eqref{spaces_cprob_LW}.
\end{proposition}
{\it Proof.}
Assuming 
\begin{equation}\label{growthphi_LW}
|h(x)-h(y)|\leq C|x-y|^{1-\epsilon}(1+|x|^{4/3}+|y|^{4/3}),\forall x,y\in\R
\end{equation}
for some $C>0, 0<\epsilon<1$, we have
\begin{align*}
|\langle h(u_n)-h(u),v\rangle_{\Wc,\Wc^*}|\leq \underbrace{C(\|u\|^\frac{4}{3}_{C(L^2)}, \|u_n\|^\frac{4}{3}_{C(L^2)})}_{<\infty}\underbrace{\|u_n-u\|^{1-\epsilon}_{L^2(L^{6-6\epsilon})}}_{\to0}\|v\|_{L^2(L^6)} \quad\to0.
\end{align*}
Next,
\begin{align*}
&|\langle f_n(\alpha_n,u_n)-f(\alpha,u),v\rangle_{\Wc,\Wc^*}|\\
&\quad=\left|\int_0^T\int_\Omega (f_n-f)(\alpha,u)v\,dx\,dt + \int_0^T\int_\Omega (f_n(\alpha_n,u_n)-f_n(\alpha,u))v\,dx\,dt\right|\\
&\quad\leq \underbrace{\left|\int_0^T\int_\Omega (f_n-f)(\alpha,u)v\,dx\,dt\right|}_{=:A_n}+\underbrace{\|(f_n)'_{\alpha,u}\|_{L^\infty(\R^{n+1})}}_{<\infty}\underbrace{(|\alpha_n-\alpha|+\|u_n-u\|_{L^2(L^{2})})}_{\to0}\|v\|_{L^2(L^2)}.
\end{align*}
In $A_n$, for fixed $u\in\Vc$ and each $ v\in\Wc\subset L^1(\tom)$, we observe that $\mu_v\in (L^\infty(\R^{n+1}))^*$ with $\|\mu_v\|=\|v\|_{L^1(\tom)}$ by defining $\mu_v:=\int_0^T\int_\Omega (\cdot)(\alpha,u)v\,dx\,dt$. Since $f_n\overset{L^\infty(\R^{n+1})}{\rightharpoonup}f$, it yields $A_n=\mu_v(f_n-f)\to0.$ Now, the bilinear term is estimated as
\begin{align}\label{Land-dis-c}
&|\langle c_nu_n-cu,v\rangle_{\Wc,\Wc^*}|=\left|\int_0^T\int_\Omega c_n(u_n-u)v\,dx\,dt+\int_0^T\int_\Omega (c_n-c)uv\,dx\,dt \right| \nonumber\\
 &\quad\leq \underbrace{\|c_n\|_{L^2}}_{<\infty}\underbrace{\|u-u_n\|_{L^2(L^3)}}_{\to0}\|v\|_{L^2(L^6)} + \Big|\int_\Omega \underbrace{(c_n-c)}_{ \rightharpoonup\, 0 \text{ in } L^2(\Omega)}\underbrace{\int_0^T uv\,dt}_{\in L^2(\Omega)}\,dx \Big| \quad\to0.
\end{align}
Weak continuity of the remaining part $(u,\source,u_0)\mapsto (\dot{u}-\Delta u-\source,u(0)-u_0)$ is straightforward, as it is a linear, bounded operator from $\Vc\times X_\source\times H$ to $\Wc\times H.$ Altogether, we claim weak continuity of $\x\mapsto (\Fb(\x)-\y):=\wtil\in \Wc,$ thus of $\x\mapsto I_V\wtil\in L^2(0,T;V).$ 
\begin{flushright}
$\diamondsuit$
\end{flushright}
With this result, we now study the week sequential continuity of $\x\mapsto-\Fb'(\x)^*(\Fb(\x)-\y)$.
\paragraph*{Weak continuity of $\x\mapsto\tilde{\mu}$ in \eqref{Land-adjoint}.}
For $\lambda=(c,\source)$, $\mutil=(\mutil_c,\mutil_\source)$, $\mutil_c=-\int_0^T F_c(\lambda,u)^*I_V\wtil\, dt$ where $\wtil=\Fb(\x)-\y$ as above, we write
\begin{align*}
&|\langle\mutil_c^n-\mutil_c,\ctil\rangle_{X_c}|:=\left|\int_0^T\int_\Omega (u_n I_V \wtil^n-u I_V\wtil)\ctil\,dx\,dt\right|\\
&\quad\leq \underbrace{\|u_n-u\|_{L^2(L^{6-\epsilon})}}_{\to0}\underbrace{\|(I_V\wtil^n)\ctil\|_{L^2(L^{(6-\epsilon)/(5-\epsilon)})}}_{<\infty} + \left|\int_0^T\int_\Omega \underbrace{u\ctil}_{\in L^2(V^*)} \underbrace{I_V (\wtil^n-\wtil)\,dx\,dt}_{\rightharpoonup\, 0 \text{ in } L^2(V)}\right| \quad\to0
\end{align*}
for any $\ctil\in X_c$,
thus showing weak continuity of $\x\mapsto\mutil_c(\x).$ Weak continuity of $\x\mapsto\mutil_\source(\x)$ could be obtained in a similar way, replacing $F_\source(\lambda,u)^*=\text{Id}$.

\paragraph*{Weak continuity of $\x\mapsto\tilde{v_0}$ in \eqref{Land-adjoint}.}
As $(u,u_0)\mapsto \tilde{v_0}:=u(0)-u_0$ is linear and bounded, its weak continuity is clear.

\paragraph*{Weak continuity of $\x\mapsto\betil$ in \eqref{Land-adjoint}.}
We consider
\begin{align}\label{Land-dis-f}
&|\langle \betil^n-\betil,\zeta\rangle_{R^n}|:=\Big|-\int_0^T\int_\Omega (f_n)'_\alpha(\alpha_n,u_n)\,\zeta\,I_V\wtil^n\,dx\,dt+\int_0^T\int_\Omega f'_\alpha(\alpha,u)\,\zeta\, I_V\wtil\,tx\,dt \Big| \nonumber\\
&=\Big|-\int_0^T\int_\Omega [(f_n)'_\alpha(\alpha_n,u_n)-(f_n)'_\alpha(\alpha,u)]\,\zeta\,I_V\wtil^n\,dx\,dt - \int_0^T\int_\Omega f'_\alpha(\alpha,u)\,\zeta\,I_V(\wtil^n-\wtil)\,dx\,dt \nonumber\\&\qquad -\int_0^T\int_\Omega [(f_n)'_\alpha(\alpha,u)-f'_\alpha(\alpha,u)]\,\zeta\,I_V\wtil^n\,dx\,dt \Big| \nonumber\\
&\leq \underbrace{\|(f_n)'_\alpha\|_{W^{1,\infty}(\R^{n+1})}}_{<\infty}\underbrace{(|\alpha-\alpha_n|+\|u-u_n\|_{C(L^2)})}_{\to0}\underbrace{\|\zeta I_V\wtil^n\|_{L^1(L^2)}}_{<\infty}\\&+ \left|\int_0^T\int_\Omega \underbrace{f'_\alpha(\alpha,u)\,\zeta}_{\in L^2(V^*)}\underbrace{I_V(\wtil^n-\wtil)}_{\rightharpoonup\,0 \text{ in } L^2(V)}\,dx\,dt\right|\\ 
&\hspace*{3cm}+\underbrace{\|\zeta I_V\wtil^n\|_{L^2(L^2)}}_{<\infty}\underbrace{\sqrt{\int_0^T\int_\Omega |(f_n)'_\alpha(\alpha,u)-f'_\alpha(\alpha,u)|^2\,dx\,dt}}_{:=A_n'}. \nonumber
\end{align}
Regarding $A_n'$, let us fix $(\alpha,u)$, then  set $\Omega_T:=\{(t,x)\in\tom: |u(t,x)|<\infty\}$. Now $\Omega_T$ has nonzero measure, as $u\in\Vc\subset L^1(\tom)$. Moreover, $|(\tom)\backslash\Omega_T|=0$. Next, for each $(t,x)\in\Omega_T$, the functional defined by $\mu_{t,x}:= (\cdot)(\alpha,u(t,x))$ belongs to $C_b(\R^{n+1})^*$ with $\|\mu_{x,t}\|=1$. From this, we ascertain \[ (f_n-f)'_\alpha(\alpha,u)(x,t)=(f_n-f)'_\alpha(\alpha,u(x,t))=\mu_{t,x}((f_n-f)'_\alpha)\to 0 \quad \forall^\text{a.e. } (t,x)\in\tom\] for $(f_n)'_\alpha\rightharpoonup f'_\alpha$ in $C_b(\R^{n+1})$. In addition, \[\sup_n\|(f_n)'_\alpha(\alpha,u)\|_{L^\infty(\tom)}+\|f'_\alpha(\alpha,u)\|_{L^\infty(\tom)}\leq \sup_n\|(f_n)'_\alpha\|_{L^\infty(\R^{n+1})}+\|f'_\alpha\|_{L^\infty(\R^{n+1})}<\infty.\]
Applying the Dominated Convergence Theorem yields $A_n'=\|(f_n)'_\alpha(\alpha,u)-f'_\alpha(\alpha,u)\|_{L^2(\tom)}\to0$. Note that this argument remains valid even for $\|\cdot\|_{L^p(\tom)}, 1\leq p<\infty$. This demonstrates the weak convergence of $\x\mapsto\betil(\x)$.

\paragraph*{Weak continuity of $\x\mapsto \tilde{v}=I_V^{-1}\ztil$ in \eqref{Land-adjoint}.}
In the first step, testing \eqref{Land-adjoint-wave} with $\Delta\ztil$ where $\ztil\in \Vc':=L^2(0,T;H^3(\Omega)\cap H^2_0(\Omega))\cap H^1(0,T;H^1(\Omega))$, yields
\begin{align}\label{Land-dis-energy}
&\|\nabla \dot{\ztil}\|_{L^2(L^2)}^2+\|\nabla\Delta\ztil\|_{L^2(L^2)}^2 \nonumber\\
&\qquad=\langle\srcWave,\Delta z\rangle_{\Vc^*,\Vc}+\int_\Omega \left(I_V\wtil(t)-\dot{\wtil}(t)\right)\Delta z(t)|^{t=T}_{t=0}\,dx - \int_0^T\int_\Omega I_V\wtil\Delta\dot{\ztil}\,dx\,dt \nonumber\\
&\qquad\leq \|\srcWave\|_{\Vc^*}\|\Delta\ztil\|_\Vc+\|\htil\|_{L^2}\|\Delta \ztil(0)\|_{L^2} +\|\nabla I_V\wtil\|_{L^2(L^2)}\|\nabla\dot{\ztil}\|_{L^2(L^2)} \nonumber\\
&\qquad\leq\left(\|\srcWave\|_{\Vc^*}+C_{\Vc'\to C(H^2)}\|\htil\|_{L^2} + \|I_V\wtil\|_{L^2(V)} \right)\|\ztil\|_\Vc',
\end{align}
where $\srcWave:=F'_u(\lambda,u)^*I_V\wtil +f_u(\alpha,u)^*I_V\wtil -M^*\ytil$ is the right hand side of the wave equation \eqref{Land-adjoint-wave} without the first term $I_V\dot{\wtil}$. Here $F'_u(\lambda,u)=-\Delta + c+h'(u)$ under the assumption $|h'(x)-h'(y)|\leq C|x-y|^{1-\epsilon}(1+|x|^{1/3}+|y|^{1/3}),\forall x,y\in\R$.

As previously, when $\x_n\rightharpoonup \x$, one has $I_V\wtil^n\rightharpoonup I_V\wtil$ in $L^2(0,T;V)$ and $\htil_n\rightharpoonup\htil$ in $L^2(\Omega)$. We now show $\srcWave_n \rightharpoonup \srcWave$ in $\Vc^*$. Indeed,
\begin{align*}
&\langle -\Delta^*(I_V\wtil^n-I_V\wtil,v\rangle_{\Vc^*,\Vc}=-\int_0^T\int_\Omega \underbrace{I_V(\wtil^n-\wtil)}_{\rightharpoonup\, 0 \text{ in } L^2(0,T;V)}\underbrace{\Delta v}_{\in\Wc}\,dx\,dt \to 0,\\[1ex]
&\langle c_nI_V\wtil^n-cI_V\wtil,v\rangle_{\Vc^*,\Vc}\to 0 \qquad \text{similarly to } \eqref{Land-dis-c} \text{ with }  I_V\wtil \text{ in place of } u,\\[1ex]
&|\langle  h'(u_n)^* I_V\wtil^n-h'(u)^* I_V\wtil),v\rangle_{\Vc^*,\Vc}| \\
&\quad\leq \underbrace{\| h'(u_n)-h'(u)\|_{L^2(L^3)}}_{\leq C(\|u\|^{\frac{1}{3}}_{C(L^2)})\|u_n-u\|^{1-\epsilon}_{L^2(L^{6-6\epsilon})}\to 0}
\underbrace{\|I_V\wtil^n\|_{L^2(L^6)} \|v\|_{C(L^2)}}_{<\infty}+\left|\int_0^T\int_\Omega \underbrace{I_V(\wtil^n-\wtil)}_{\rightharpoonup\, 0 \text{ in } L^2(0,T;V)} \underbrace{h'(u)v}_{\in\Wc}\,dx\,dx\right|,\\[1ex]
&\langle f'_u(\alpha_n,u_n)^*I_V\wtil^n-f'_u(\alpha,u)^*I_V\wtil,v\rangle_{\Vc^*,\Vc}\to 0\\ &\qquad\text{similarly to } \eqref{Land-dis-f} \text{ with } f'_u \text{ in place of } f'_\alpha, v \text{ in place of } \zeta.
\end{align*}
The last estimate is analogous to \eqref{Land-dis-f}, but modifies  the upper bound for the term involving $A_n'$ to $\|v\|_{C(L^2)}\|I_V\wtil^n\|_{L^2(L^6)}\|(f_n)'_u(\alpha,u)-f'_u(\alpha,u)\|_{L^2(L^3)}$. This yields $\srcWave_n \rightharpoonup \srcWave$ in $\Vc^*$ when $\x\rightharpoonup 0$, as claimed.

In order to form the full $\Vc'$-norm on the left hand side of \eqref{Land-dis-energy}, we test \eqref{Land-adjoint-wave} by $\ztil$. By then applying Young's inequality with $\epsilon>0$, we eventually obtain
\begin{align}\label{Land-dis-energy1}
(1-3\epsilon)\|\ztil_n-\ztil\|^2_{\Vc'}\leq \frac{1}{\epsilon}\left(\|\srcWave_n-\srcWave\|_{\Vc^*}^2+(C_{\Vc'\to C(H^2)})^2\|\htil_n-\htil\|^2_{L^2} + \|I_V(\wtil^n-\wtil)\|^2_{L^2(V)} \right).
\end{align}
Using Galerkin approximation, one can show that for each $(\srcWave,\htil,\wtil)\in \Vc^*\times L^2(\Omega)\times \Wc$, there exists a unique $\ztil\in \Vc'$ solving \eqref{Land-adjoint-wave}. Moreover, $\ztil$ depends continuously on the data $(\srcWave,\htil,\wtil)$ through the expression \eqref{Land-dis-energy1}. Since $\Vc^*\times L^2(\Omega)\times \Wc\ni(\srcWave,\htil,\wtil)\mapsto\ztil\in \Vc'$ is linear and bounded, it is weakly continuous. In conclusion, when $\x_n\rightharpoonup \x$, we have $\ztil_n {\rightharpoonup} \ztil$ in $\Vc'$, equivalently $I_V^{-1}\ztil_n {\rightharpoonup} I_V^{-1}\ztil$ in $\Vc$, proving weak continuity of $\bfx\mapsto \vtil:=I_V^{-1}\ztil$ in \eqref{Land-adjoint}.

\paragraph*{Weak continuity of $\x\mapsto\gtil$ in \eqref{Land-adjoint}.}
For $\gtil$ as in \eqref{Land-adjoint} and setting $C_\pi:=\frac{1}{2\pi^{(n+1)/2}}$, we evaluate, for any $\psi\in\Cc$
\begin{align*}
&\langle \gtil^n-\gtil,\psi\rangle_\Cc=\int_{\mathbb{R}^{n+1}} (1+|\kappa|^2)^\ell
(1+|\omega|^2)^r \overline{\Fc\psi(\kappa,\omega)}\Fc(\gtil^n-\gtil)(\kappa,\omega)\, d\kappa\, d\omega\\
&=C_\pi\int_{\mathbb{R}^{n+1}}\overline{\Fc\psi(\kappa,\omega)}\left[\int_0^T\int_\Omega e^{-i\kappa\cdot\alpha_n-i\omega u_n(x,t)}I_V\wtil^n-e^{-i\kappa\cdot\alpha-i\omega u(x,t)}I_V\wtil\,dx\,dt \right]\,d\kappa\,d\omega \\
&=C_\pi\int_{\mathbb{R}^{n+1}}\overline{\Fc\psi(\kappa,\omega)}\Bigg[-\int_0^T\int_\Omega \underbrace{e^{-i\kappa\cdot\alpha-i\omega u}}_{\in L^2(V^*)}\underbrace{I_V(\wtil^n-\wtil)}_{\rightharpoonup\,0 \text{ in } L^2(V)}\,dx\,dt \\
&\quad+\int_0^T\int_\Omega \underbrace{\int_0^1 e^{-i\kappa\cdot(\alpha-\theta(\alpha_n-\alpha))-i\omega (u+\theta(u_n-u))}\,d\theta}_{\in[-1,1]} \underbrace{(-i\kappa\cdot(\alpha_n-\alpha)-i\omega(u_n-u))I_V\wtil^n}_{\to\, 0 \text{ in } L^1(0,T;L^1(\Omega))}\,dx\,dt \Bigg]\,d\kappa\,d\omega\\
&=:C_\pi\int_{\mathbb{R}^{n+1}}\overline{\Fc\psi(\kappa,\omega)} B_n(\kappa,\omega)\,d\kappa\,d\omega.
\end{align*}
and deduce pointwise convergence in $(\kappa,\omega)$ of $B_n$. Together with uniform boundedness via
\begin{align*}
&\int_{\mathbb{R}^{n+1}}\overline{\Fc\psi(\kappa,\omega)} B_n(\kappa,\omega)\,d\kappa\,d\omega\\ 
&\qquad\leq \int_{\mathbb{R}^{n+1}}|(1+|\kappa|^2)^{1/2}(1+|\omega|^2)^{1/2}\overline{\Fc\psi(\kappa,\omega)}|\underbrace{\left(\frac{C}{(1+|\kappa|^2)^{1/2}(1+|\omega|^2)^{1/2}}\right)}_{\in L^2(\R^{n+1})}\,d\kappa\,d\omega\\
&\qquad\leq C\|\psi\|_{H^1(\R^{n+1})}\leq C\|\psi\|_\Cc,
\end{align*} 
and application of the Dominated Convergence Theorem, we conclude $\langle \gtil^n-\gtil,\psi\rangle_\Cc\to0$, yielding weak continuity of $\x\mapsto\gtil(\x)$.

\subsubsection{Lipschitz continuity of $\bfx\mapsto\Fb'(\bfx)^*(\Fb(\bfx)-\bfz)$}\label{sec:Lipschitz}
Above, we have verified weak continuity of $\bfx\mapsto-\Fb'(\bfx)^*(\Fb(\bfx)-\y)$, where all the estimates were written in the form of $\|\bfx_n-\bfx\|$. Therefore, Lipschitz continuity of $\Fb'(\bfx)^*(\Fb(\bfx)-\bfz)$ could be established in the same manner, the Lipschitz constant $L<2$ being obtained in the ball $B_R(\bfx^\dagger)$ with sufficiently small $R.$

\medskip

As a consequence, we can conclude from Propositions~\ref{prop:weakcont}, \ref{prop:LW} the following convergence results on Tikhonov regularization and Landweber iteration.
\begin{corollary}[Tikhonov]\label{convTikh_1}
For the operator $\Fb$ defined by \eqref{Land-dis-eq} on the spaces \eqref{VWY_LW}, \eqref{spaces_cprob_LW}, \eqref{spaces_cprob_LW_1}, $\Cc=H^s(\mathbb{R})$, 
\blue{$s>\frac52$} with $M:\Vc\to\Yc$ linear and bounded and $h$ satisfying \eqref{growthphi_LW} we have subsequential convergence of 
$(c^{\gamma(\delta),\delta,N(\delta)},\source^{\gamma(\delta),\delta,N(\delta)},u_0^{\gamma(\delta),\delta,N(\delta)},u^{\gamma(\delta),\delta,N(\delta)}, f^{\gamma(\delta),\delta,N(\delta)})$ to a solution of the inverse problem \eqref{IP_c-problem}, \eqref{Muy}
as $\delta\to0$.
\end{corollary}

\begin{corollary}[Landweber]
Let the assumptions of Corollary~\ref{convTikh_1} hold and additionally assume full measurements $M=\mbox{id}_{\Vc\to\Yc}$ and \eqref{klogdelta}. Then we have weak subsequential convergence of
$(c_{N(\delta),k_*(\delta)}^\delta,\source_{N(\delta),k_*(\delta)}^\delta,{u_0}_{N(\delta),k_*(\delta)}^\delta,u_{N(\delta),k_*(\delta)}^\delta, f_{N(\delta),k_*(\delta)}^\delta)$ to a solution of the inverse problem \eqref{IP_c-problem}, \eqref{Muy} as $\delta\to0$.
\end{corollary}

\section{Outlook}
In this study, we have carried out a convergence analysis for discretizations of Tikhonov and projected Landweber regularization using Neural Networks. The convergence analysis is based on \blue{a} priori choice of the regularization and discretization parameters, \blue{where the latter relates to the network approximation error. Our analysis is applicable not only for discretization by NNs, but also for general discretization schemes.} As an application, we have presented a parameter identification problem for a time-dependent PDE, whose unknown nonlinearity is approximated by a neural network.
Our all-at-once approach does not require a training process for learning the nonlinearities beforehand, instead simultaneously determining it alongside the unknown coefficients and the solution of the PDE.

This paper focuses on the theoretical aspects. Numerical results for the regularization with neural networks can be found in \cite{Nguyen:21}. \blue{Also in \cite{Nguyen:21}, further details on the discretized problem are discussed, such as differentiability of the forward mapping, unique existence for the \emph{learning-informed PDEs} (NN as a reaction term in the PDE), the tangential cone condition for networks and so forth. A potential extension to our study is the inclusion of further components in the unknown nonlinear response, e.g. $f(\alpha,u,\nabla u,\nabla^2 u\ldots)$, for more flexible models, as was done in \cite{ Brunton16,Long18,Raissi18}.} 

On the analytical side, an open problem is determining convergence rates for Tikhonov regularization, based on quantified approximation results for neural networks \blue{as in \cite{BurgerEngl20} (see also Discussion \ref{appox-rate} in the introduction)}. These rates require enhanced regularity of the exact solution in terms of so-called source conditions, whose interpretation for the problem setting considered here is another interesting task.

\paragraph{Acknowledgments.}
The work of the first author was supported by the Austrian Science Fund {\sc fwf} under the grants P30054 and DOC 78. Moreover, we wish to thank both reviewers for fruitful comments leading to an improved version of the manuscript.

\bibliography{lit}{}
\bibliographystyle{plain} 
\end{document}